\documentclass[12pt,a4paper]{amsart}
\ifx\undefined\pdfpagewidth
\usepackage{hyperref}
\else
\usepackage[pdftex]{hyperref}
\fi
\hypersetup{%
  colorlinks,
  bookmarksopen,
  bookmarksnumbered,
  citecolor=blue,
  linkcolor=black,
  pdfstartview=FitH,
}
\usepackage{amssymb, amstext, amscd, amsmath, color}
\usepackage{graphicx}
\usepackage{url}
\usepackage{multirow}
\theoremstyle{plain}
\newtheorem{thm}{Theorem}[section]
\newtheorem{cor}[thm]{Corollary}
\newtheorem{prop}[thm]{Proposition}
\newtheorem{lem}[thm]{Lemma}

\theoremstyle{definition}
\newtheorem{rem}[thm]{Remark}
\newtheorem{defn}[thm]{Definition}

\newtheorem{eg}[thm]{Example}

\renewcommand{\phi}{\varphi}

\newcommand{\ol}{\overline}
\newcommand{\ul}{\underline}

\newcommand{\re}{\operatorname{Re}}

\newcommand{\upchi}{{\raise.35ex\hbox{$\chi$}}}

\newcommand{\qand}{\q{and}}
\newcommand{\qor}{\q{or}}

\newcommand{\q}[1]{\quad\text{#1}\quad}
\newcommand{\Alg}{\operatorname{Alg}}

\newcommand{\diag}{\operatorname{diag}}

\newcommand{\Lat}{\operatorname{Lat}}

\newcommand{\Unit}{\operatorname{Unit}}
\newcommand{\bC}{{\mathbb{C}}}
\newcommand{\bD}{{\mathbb{D}}}
\newcommand{\bN}{{\mathbb{N}}}

\newcommand{\bR}{{\mathbb{R}}}
\newcommand{\bT}{{\mathbb{T}}}

\newcommand{\bZ}{{\mathbb{Z}}}
\newcommand{\bH}{{\mathbb{H}}}
  \newcommand{\A}{{\mathcal{A}}}
  \newcommand{\B}{{\mathcal{B}}}
  
  \newcommand{\D}{{\mathcal{D}}}

\renewcommand{\H}{{\mathcal{H}}}
  \newcommand{\I}{{\mathcal{I}}}
  
  \newcommand{\K}{{\mathcal{K}}}  
\renewcommand{\L}{{\mathcal{L}}}
  \newcommand{\M}{{\mathcal{M}}}
  \newcommand{\N}{{\mathcal{N}}}

\renewcommand{\P}{{\mathcal{P}}}

\renewcommand{\S}{{\mathcal{S}}}
  
  \newcommand{\U}{{\mathcal{U}}}
  
  \newcommand{\W}{{\mathcal{W}}}

\newcommand{\equator}{\mathrm{e}}
\newcommand{\Sigmae}{\Sigma_\equator}
\newcommand{\Sigmaa}{\Sigma_\analytic}
\newcommand{\Sigmav}{\Sigma_\volterra}
\newcommand{\FPS}{\Sigma_{\textsc{FP}}}

\newcommand{\analytic}{\mathrm{a}}
\newcommand{\volterra}{\mathrm{v}}
\newcommand{\volterranest}{\N_\volterra}
\let\vnest\volterranest

\newcommand{\analnest}{\N_\analytic}
\let\anest\analnest

\newcommand{\Nv}{\volterranest}

\newcommand{\nests}{\mathfrak{N}}
\newcommand{\Ltp}{\Ltwo{\bR_+}}
\newcommand{\Ltwo}[1]{\Leb{2}{#1}}
\newcommand{\Lt}{\Ltwo{\bR}}
\newcommand{\Ltm}{\Ltwo{\bR_-}}
\newcommand{\Htwo}[1]{\Hp{2}{#1}}
\newcommand{\Ht}{\Htwo{\bR}}
\newcommand{\Leb}[2]{L^{#1}(#2)}
\newcommand{\Hp}[2]{H^{#1}(#2)}
\newcommand{\Htb}{\overline{\Ht}}
\newcommand{\sot}{strong operator topology}
\newcommand{\Linf}[1]{\Leb{\infty}{#1}}
\newcommand{\Li}{\Linf{\bR}}
\newcommand{\Ad}{\operatorname{Ad}}
\newcommand{\supp}{\operatorname{supp}}
\newcommand{\Proj}{\operatorname{Proj}}

\newcommand{\sgn}{\operatorname{sgn}}
\newcommand{\bB}{\mathbb{B}}

\DeclareMathOperator*{\sotlim}{{\text{\sc sot}}-lim}
\newcommand{\Mhyp}{\hat\L}
\newcommand{\Uhyp}{\U(\Mhyp)}

\newcommand{\hyp}{{\mathrm{hyp}}}

\newcommand{\Shyp}{\Sigma_\hyp}
\newcommand{\bdy}{\partial\M}
\newcommand{\clo}[1]{\overline{#1}}

\newcommand{\sotto}{\stackrel{\textsc{sot}}{\to}}
\let\lt\Lt
\let\oldepsilon\epsilon
\let\epsilon\varepsilon
\let\varepsilon\oldepsilon
\newcommand{\sd}{\mathbin{\text{\footnotesize$\bigtriangleup$}}}
\newcommand{\UHP}{\bH}
\newcommand{\bd}{\partial}
\newcommand{\pure}{\textrm{pure}}
\newcommand{\ub}[1]{^{(#1)}}
\newcommand{\bnc}{\clo{\bB^n}}
\newcommand{\sn}{S^{n-1}}
\newcommand{\bmodsim}[1][\bnc]{#1/\!\!\sim}
\newcommand{\cd}{\clo\bD}
\newcommand{\phis}[1][s]{\phi_{#1}}
\newcommand{\Mphis}[1][s]{M_{\phis[#1]}}
\newcommand{\bfn}[1]{#1^{\star}}
\begin{document}
\def\nl{\par}
\def\keywords{Projection manifold, Beurling subspace manifold, locally
  unitary manifold, Hardy space, strange
  limit, Fourier-Plancherel sphere, hyperbolic sphere}
\thanks{%
  \today.\nl
  {\it $2000$ Mathematics Subject Classification}. 
  47B38,            % operators on function spaces
  46E20,            % Hilbert spaces of continuous, 
                    %     differentiable or analytic functions
                    %
  58D15.            %  analysis on manifolds -> manifolds of mappings
  \nl
  }
\title{Manifolds of Hilbert space projections}
\hypersetup{pdftitle={Manifolds of Hilbert space projections}}
\author[R.H.Levene]{R. H. Levene}
\address{Queen's University Belfast\\ Belfast\; BT7~1NN\\UK}
\email{r.levene@qub.ac.uk}
\author[S.C.Power]{S. C. Power}
\hypersetup{pdfauthor={R.H. Levene and S.C. Power}}
\hypersetup{pdfkeywords={\keywords}}
\hypersetup{pdfsubject={Functional analysis}}
\address{Lancaster University\\ Lancaster\; LA1~4YF\\UK}
\email{s.power@lancaster.ac.uk}

\begin{abstract}
  The Hardy space $\Ht$ for the upper half plane together with a
  unimodular function group representation $u(\lambda) =
  \exp(i(\lambda_1\psi_1 + \dots + \lambda_n\psi_n))$, $\lambda \in
  \bR^n$, gives rise to a manifold $\M$ of orthogonal projections for
  the subspaces $u(\lambda)\Ht$ of $\Lt$. For classes of
  admissible functions $\psi_i$ the strong operator topology closures
  of $\M$ and $\M \cup \M^\bot$ are determined explicitly as various
  $n$-balls and $n$-spheres. The arguments used are direct and rely on
  the analysis of oscillatory
  integrals (Stein~\cite{stein93:harmonic-analysis}) and Hilbert space
  geometry. Some classes of these closed projection manifolds are classified
  up to unitary equivalence. In particular the Fourier-Plancherel
  $2$-sphere and the hyperbolic $3$-sphere of Katavolos and
  Power~\cite{pow-kat:hyp} appear as distinguished special cases
  admitting nontrivial unitary automorphisms groups which are explicitly
  described.
\end{abstract}

\maketitle

\section{Introduction}

Let $\M$ be a set of closed subspaces of a Hilbert space $\H$
endowed with the strong operator topology inherited from the
identification of closed subspaces~$\K$ with their self-adjoint
projections $[\K]:\H\to \K$. If $\M$ is finitely parametrised in
the sense that
\[
\M = \{[\K_\lambda] : \lambda \in M \subseteq \bR^n\}
\]
with $M$ a topological manifold, then $\M$ may in fact be
homeomorphic to $M$. Furthermore $\M$ may admit a certain local
unitary description and an associated smooth structure under which
it is a diffeomorph of a differentiable manifold in $\bR^n$.
Natural examples of such manifolds of projections are provided by
Grassmannian manifolds and their submanifolds. Also, operators $T$
in the first Cowen-Douglas class~\cite{cowen-douglas},
\cite{curto-salinas} for a complex connected domain $\Omega$ in
$\bC^m$ provide diverse realisations (even for $m=1$) of domains
in $ \bR^{2m}$, namely
\[
\M_T = \{[\ker (T-\omega I)]: \omega \in \Omega \subseteq
\bR^{2m}\}.
\]
However our primary motivating examples derive from invariant
subspaces for semigroups of unitary operators, such as the
semigroup $\W=\{\alpha U_tV_s : s,t \in \bR_+,\ |\alpha| = 1\}$
associated with jointly irreducible one parameter unitary groups
satisfying the Weyl commutation relations. Furthermore, such
subspaces and their complementary spaces are generally infinite
dimensional.  These examples motivate the consideration of general
subspace manifolds as formulated in
Definitions~\ref{def:Cinfsubman}, \ref{def:locunitsubman}
and~\ref{def:beurlingsubman}.

The embracing realm we consider is the traditional one allied to
operator function theory, namely the set of subspaces of $\Lt$ of
the form
\[
u\Ht, \quad u\Htb, \quad \Ltwo{E}
\]
where $\Ht$ is the Hardy space for the upper half plane, $u(x)$ is
unimodular, and  $L^2(E)$
 is the  space of functions supported
on a measurable set $E$. We examine subspace manifolds of the form
\[ \M = \M(\S) = \{e^{i\psi(x)}\Ht : \psi(x) \in \S\},
\]
where $\S$ is a finite dimensional real vector space of
real-valued functions, we analyse limits of projections and we
identify the associated closed topological manifolds. Our approach
and the ensuing identifications give a unified explanation for
various so called ``strange limits'' of projections. These include
the special cases considered by Katavolos and
Power~\cite{pow-kat:FB}, \cite{pow-kat:hyp} which were derived by ad hoc
arguments leaning on operator algebra methods. Specifically we
show by direct methods that the space of functions
\[
\S_1= \{\lambda_1x +\lambda_2x^2 : \lambda = (\lambda_1,\lambda_2)
\in \bR^2\},
\]
has  subspace manifold $\M(\S_1)$ whose closure is homeomorphic to
the closed unit disc, while for the space
\[
\S_2= \{\lambda_1\log |x|  +\lambda_2x + \lambda_3x^{-1} : \lambda
 \in \bR^3\},
\]
the manifold $\M(\S_2)$ has closure homeomorphic to the closed unit ball in~$\bR^3$.
In contrast, the closures of Cowen-Douglas projection
manifolds are generally one point compactifications.

A consequence of the limit projection analyses
in~\cite{pow-kat:FB}, \cite{pow-kat:hyp} is that a reducing invariant
subspace for the Weyl semigroup $\W$, or for the $ax+b$ unitary
semigroup (with $a\ge 1$ and $b\ge 0)$, turns out to be a strong
operator topology limit of a sequence of purely invariant
projections, that is, a limit of those with no reducing part. We
obtain in Theorem~\ref{thm:purelimits} a similar phenomenon for
the multiplication semigroup $\{M_{e^{i\lambda x}}:\lambda \ge
0\}$ acting on $L^2(\bR)$. Equivalently, translating to the circle
$\bT$, we show that each reducing invariant subspace for the
bilateral shift, which corresponds to a measurable subset of
$\bT$, is a strong operator topology limit of Beurling projections
$[uH^2(\bT)]$. We are not aware of any other proof, direct or
indirect, of this seemingly classical fact.

\begin{figure}
  \label{fig:FPS}
  \centering
  \includegraphics{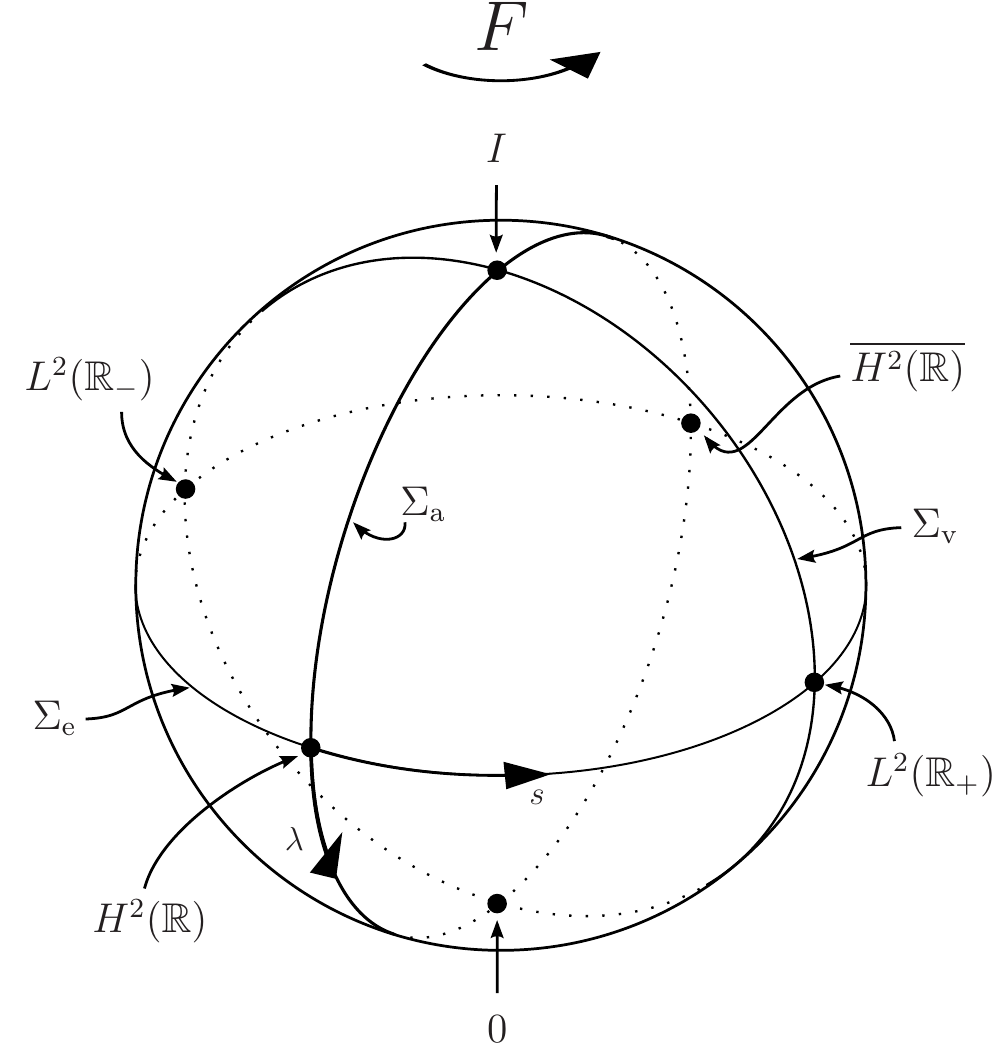}
  \caption{A natural realisation of $\FPS$, the Fourier-Plancherel
    sphere, on which the Fourier-Plancherel transform $F$ acts as a
    quarter-rotation.}
\end{figure}

Our principal tool is Theorem~\ref{thm:psi} which, for a
function~$\psi$ in a certain
admissible class, identifies the limit of the
projections $[e^{in\psi}\Ht]$, as $n \to \infty$,
 with the projection $[L^2((\psi')^{-1}(-\infty,0))]$. The proof
makes use of methods of Stein~\cite{stein93:harmonic-analysis} for
the analysis of oscillatory integrals. We go on to show that for
quite general $n$ dimensional spaces $\S$ of admissible functions
the subspace manifold $\M(\S)$ has closure homeomorphic to the
closed unit ball in $\bR^n$. Moreover, in these cases  the
two-component subspace manifold $\M(\S) \cup \M(\S)^\perp$ has
closure, denoted  $\Sigma(\S)$, which is homeomorphic to an
$n$-sphere. The $2$-sphere $\FPS = \Sigma(\S_1)$ is the so-called
Fourier-Plancherel sphere of \cite{pow-kat:hyp}. (See
Figure~\ref{fig:FPS} and Example~\ref{eg:FPS}.) It is natural to
consider how such Hilbert space manifolds may be classified
geometrically, that is, up to unitary equivalence, and in the
examples of Section~\ref{sec:4} and Section~\ref{sec:5.3} we distinguish a number of
distinct 2-spheres and 3-spheres.

 In Section~\ref{sec:5} we consider the sphere
$\Sigma(\S_1)$ and its hyperbolic variant, the 3-sphere
$\Sigma(\S_2)$, from the point of view of their unitary
automorphism symmetries. The analysis here exploits the nontrivial
foliations induced by the natural order on projections. It is
shown that the group $\U(\Sigma(\S_1))$ of unitaries which act
bijectively on $\Sigma(\S_1)$ is generated by the set
\[
\{M_{e^{i\lambda x}}, M_{e^{isx^2}}, F, \alpha I: \lambda, s \in
\bR, \ |\alpha | = 1\},
\]
that is, by two 1-parameter unitary groups, the scalar circle
group and the Fourier-Plancherel transform $F$. In particular,
this group contains the dilation semigroup $\{V_t: t \in \bR\}$. A
similarly detailed description is given for $\U(\Sigma (\S_2))$
and this leads to the unitary automorphism group identification
as a certain double semidirect product
 \[
 \Ad(\U(\Sigma (\S_2))) \cong \big( \bR^3 \rtimes \bR\big) \rtimes
(\bZ/2\bZ)^2.
\]

The manifolds $\M(\S)$, \textit{Beurling manifolds} in our
terminology, may be regarded as smooth in the strict, locally
unitary, sense that a neighbourhood of a subspace $\K$ is given by
the local action on $\K$ of a certain unitary group representation
of $\bR^n$.
Furthermore the Fourier-Plancherel 2-sphere $\FPS$  is remarkable
in being smooth in this way at all points except the poles $0, I$.
It of interest then to identify similar compact projection
manifolds which are smooth off a finite set. In this regard we see
in Section~\ref{sec:5.3} that this is not generally the case for other polynomial
$2$-spheres $\Sigma(\S)$.

As we have intimated above our considerations lie entirely in the
realm of operator function theory tied to the Hardy space for the
line. However the oscillatory integral methods are expected to be
effective in multivariable function spaces and for higher rank
settings in noncommutative harmonic analysis. This should lead to
the identification of other closed subspace manifolds with
interesting topology and geometry.
\medskip

We wish to thank Alexandru Aleman, Gordon Blower, Jean Esterle,
Aristides Katavolos, Alfonso Montes Rodriguez and Donald Sarason for
their interest and communications concerning strange limits.

%%%%%%%%%%%%%%%%%%%%%%%%%%%%%%%%%%%%%%%%%%%%%%%%%%%%%%%%%%%
\section{Subspace manifolds}

In this section we give some definitions and examples.

Let $\H$ be a separable Hilbert space, $\Proj(\H)$ the set of
self-adjoint projections and $\Unit(\H)$ the set of unitary
operators. We shall routinely identify a closed subspace $\K$ with
its associated orthogonal projection, denoted $[\K]$.

\begin{defn}\label{def:Cinfsubman}
  (i) A \emph{topological subspace manifold} in $\B(\H)$ of dimension
  $n$ is a set $\M \subseteq \Proj(\H)$, considered with the relative
  strong operator topology, which is locally
  homeomorphic to $\bR^n$. 

  (ii) A \emph{$C^\infty$ projection manifold in $\B(\H)$} (or
  \emph{$C^\infty$ subspace manifold}) is a topological subspace
  manifold $\M$ of dimension $n$ together with an atlas of charts
  $x_i: \bR^n \to \M$ (with open domains and ranges covering $\M$) for
  which the coordinate functions $x_i^{-1}x_j$ (with nonempty domain)
  are $C^\infty$ and such that for each chart $x$ with domain $U_x$
  there is a dense subspace $\D_x$ of \emph{$C^\infty$ vectors} in
  $\H$; that is, for $f,g \in \D_x$, the function $\lambda \mapsto
  \langle x(\lambda)f,g \rangle$ is $C^\infty$ on $U_x$.
\end{defn}

In fact we bypass the technicalities of (ii) in
Definition~\ref{def:Cinfsubman} in view of the fact that the
smooth subspace manifolds we consider have a stronger locally
unitary structure as in the following formal definition.

\begin{defn}\label{def:locunitsubman}
A \emph{locally unitary subspace manifold} of dimension $n$ in
$\B(\H)$
 is a topological subspace manifold~$\M$ such that
 for each $P$ in $\M$ there is a strong operator topology
neighbourhood in $\M$ of the form
\[
\N_P = \{[\rho_P(\lambda)P\H]: \lambda \in \bB_n\}
\]
where $\rho_P : \bR^n \to \Unit(\H)$ is a strong operator topology
continuous representation which is a homeomorphism of $\N_p$ into
$\Unit(\H)$.
\end{defn}

It is well-known that such a representation $\rho_P$ does possess a
dense subspace of $C^\infty$ vectors; see~\cite{taylor}, for example.
\medskip

We shall consider Beurling subspace manifolds  $\M$, given
formally in the next definition, together with their complement
completions, by which we mean the strong operator topology closure
of $\M \cup \M^\bot$.
\medskip

\begin{defn}\label{def:beurlingsubman}
  A \emph{Beurling subspace manifold} of dimension $n$ is a
  topological subspace manifold $\M \subseteq \Proj(\Lt)$ of the
  form
  \[
  \M = \{[u_\lambda \Ht]: \lambda \in \bR^n\}
  \]
  where $\lambda \mapsto u_\lambda$ is a weak star continuous
  representation of $\bR^n$ as unimodular functions, so that $\M$ is
  locally homeomorphic to $\bR^n$ by the single chart $x(\lambda)=
  [u_\lambda \Ht]$.
\end{defn}

If $\psi$ is a non-constant real continuous function on the line
then the projections $[\exp(i\lambda\psi)\Ht]$, for $\lambda \in
\bR$, give a one dimensional topological manifold. When $\psi(x)
=x$ the closure in $\Proj(\Lt)$  adds the subspaces $\{0\}$ and
$\Lt$ and  the complement completion, $\Sigmaa$ say, is
topologically a circle. On the other hand for $\psi(x) = x^2$ we
shall see that the closure of $\M(\{x^2\})$ adds $[L^2(\bR_+)]$
and $[L^2(\bR_-)]$ and that the complement completion $\Sigmae$
is a locally unitary $ C^\infty$ subspace manifold diffeomorphic
to the circle. We see later in Example~\ref{eg:FPS} that
$\Sigmae$ and $\Sigmaa$ are, respectively, the equator and a
great circle of the Fourier-Plancherel sphere.

The subspaces $K_{\lambda,s}=e^{i\lambda x}e^{isx^2}\Ht$, for
$s<0$ and $\lambda \in \bR$, form a subspace manifold $\M$ which arises
in the analysis of the invariant subspaces for the Weyl semigroup
\[\W=\{\alpha M_\lambda D_\mu:\lambda,\mu\ge0,\ |\alpha|=1\},\]
where  $D_\mu$ is the translation unitary $D_\mu f(x)=f(x-\mu)$.
Indeed it was shown in~\cite{pow-kat:FB} that the invariant
subspaces of $\W$ are the spaces $K_{\lambda,s}$ for $s\le0$,
together with $L^2(t,\infty)$ for $t$ in
$\bR\cup\{\pm\infty\}$ and,  moreover, that the latter subspaces
are in the closure of $\M$. Extending the parameter range of $s$
to include $s\ge 0$ and taking the complement completion one obtains the
Fourier-Plancherel sphere. The \emph{Volterra circle} $\Sigmav$
consists of the subspaces $L^2(t,\infty)$, $L^2(-\infty, t)$ for
$t\in\bR\cup\{\pm \infty\}$ and is unitarily equivalent to the great circle
$\Sigmaa$ via the Fourier-Plancherel transform.

Consider now the three dimensional subspace manifold
\[
\M=\{ |x|^{is}e^{i\lambda x}e^{i\mu
x^{-1}}\Ht:(s,\lambda,\mu)\in\bR^3\}
\]
which is the Beurling manifold $\M(\S_2)$ for the space \[\S_2=\{
s\log|x| + \lambda x + \mu x^{-1} : (s,\lambda,\mu)\in\bR^3 \}.\]
In Section~\ref{sec:4} we give a new direct proof that $\clo\M$ is
homeomorphic to a closed 3-ball and that $\Sigma (\S_2) =
\clo\M\cup\clo{\M^\perp}$ is a topological 3-sphere.

In a finite dimensional Hilbert space a locally unitary manifold
can be viewed as a submanifold of a Grassmannian manifold and it is
instructive to note the following simple examples.

Let $a$ be real, let $U_t$ be the unitary operator $\diag(e^{it},
e^{iat})$ in $\B(\bC^2)$, let $\eta$ be the unit vector
$\frac{1}{\sqrt{2}}(1,1)$ and define $\P_a = \{[U_t\bC\eta]: t\in
\bR\}$. The function $t \mapsto P_t =[U_t\bC\eta]$ is periodic if
and only if $a$ is rational, in which case $\P_a$ is a locally
unitary manifold which is homeomorphic to a circle. On the other
hand if $a$ is irrational then the vectors $U_t\eta$ return
infinitely often to any neighbourhood of a particular such vector.
It follows that the relative strong operator topology does not
agree with the usual topology of  $\bR$ and that $\P_a$ is not a
locally unitary manifold.

More generally, let $\M$ be a connected locally unitary
$1$-dimensional subspace manifold on a Hilbert space $\H$, with
neighbourhood \[\N = \{[U_tP\H]: t\in (-1,1)\}\] for the projection
$P$. Assume moreover that $U_t = \exp(itA)$ with $A$ a bounded
self-adjoint operator. The derivative of the continuous
projection-valued function $P(t) = U_tPU_t^*$ is the operator
$U_tYU_t^*$ where $Y = -i(PA-AP)$. In particular, the unitary
equivalence class of the self-adjoint operator $Y$ is an invariant
for $\M$ for unitary equivalence. This underlines the fact that
unitary equivalence here is a strong form of geometric equivalence
for subspace manifolds.  In particular, using this derivative
invariant, one can soon see that the circle manifolds $\P_a, \P_b$
are unitarily equivalent if and only if the rational numbers $a,
b$ coincide.

The Beurling subspace manifolds given above might be more precisely
specified as Euclidean Beurling subspace manifolds as there are many
other interesting $C^\infty$ projection manifolds associated with
Hardy spaces and unimodular functions.  We do not develop this here
but we note some fundamental examples.

For $n =1,2,\dots $ let
\[
\M_n = \{ [u_{\ul{\lambda}}\Ht ] : \ul{\lambda} = (\lambda_1,
\dots , \lambda_n) \in U^n\}
\]
where $u_{\ul{\lambda}}(z)$ is a Blaschke factor inner function
with zeros, possibly repeated, at points $\lambda_1, \dots ,
\lambda_n$ in the upper half plane~$\UHP$. Then $\M_n$ is a
$C^\infty$ projection manifold in $\B(\Ht)$. Also, $\M_1$ consists
of codimension 1 projections and is locally unitary with respect
to a representation of the M\"obius group rather than the
Euclidean group. The closure of $\M_1$ adds one extra projection,
namely $[\Ht]$, and is a topological
 subspace manifold homeomorphic to the $2$-sphere, realised as the one point
compactification of $\UHP$.

Manifolds analogous to these, with finite or cofinite dimensional
spaces, may be defined also for Bergman Hilbert spaces, and, more
generally, in the setting of Hermitian holomorphic vector bundles
associated with operators in the Cowen-Douglas
theory~\cite{cowen-douglas}. For example, with weighted Bergman
Hilbert spaces in place of $\Ht$ one obtains projection manifold
realisations of the unit disc with one point compactification
closures. That these are unitarily inequivalent was essentially
shown in~\cite{cowen-douglas} by the construction of curvature
invariants for Hermitian holomorphic vector bundles. This reveals
once again the geometric strictness of unitary equivalence. We
note the following alternative curvature free approach to this and
in subsequent sections we find, similarly, that we do not need to
consider curvature. However, it goes without saying that it would
be interesting and useful to define curvature invariants for
general projection manifolds.

Let $\alpha \ge 0$ and let $A^2_\alpha$ be the weighted Bergman
Hilbert space of holomorphic functions in the unit disc that are
square integrable with respect to $(1-|z|^2)^\alpha dA$, where
$dA$ is area measure. For $\lambda $ in $\bD$ let $u_\lambda(z)$
be the inner function $(\lambda - z)/(1-\ol{\lambda}z)$ and let
$\M^{\alpha}$ be the projection manifold
\[
\M^{\alpha} = \{[u_\lambda A^2_\alpha] : \lambda \in \bD\}.
\]
The range of the complementary projection $[u_\lambda
A^2_\alpha]^\bot$ is one-dimensional and is spanned
by the function $(1-\ol{\lambda}z)^{2 -\alpha}$. (See Zhu~\cite{Zhu},
for example.) A standard argument with eigenvectors (see Theorem 3.6
of Thomson~\cite{thomson} for example) shows that an operator which
leaves invariant all these subspaces is necessarily multiplication by
the complex conjugate of an $H^\infty$ function.  Thus if $U$ is a
unitary with $\Ad U(\M^{\alpha}) = \M^{\beta}$ then we may deduce that
$\Ad U$ gives an isomorphism between the respective $H^\infty$
multiplication algebras. Composing with a unitary automorphism of
$\B(A^2_\beta)$ which effects the inverse M\"obius automorphism on the
range algebra we thus obtain a unitary operator $W$ between the
Bergman spaces which intertwines multiplication: $M_zW = WM_z$. Thus
if $W1 = g$ then $Wz^m = z^mg$. For $\alpha \ne \beta$ this is
contrary to $W$ being isometric.

\begin{prop}
For $\alpha, \beta \ge 0 $ the projection manifolds $\M^{\alpha}$
and $ \M^{\beta}$ are unitarily equivalent if and only if $\alpha
= \beta$.
\end{prop}

%%%%%%%%%%%%%%%%%%%%%%%%%%%%%%%%%%%%%%%%%%%%%%%%%%%%%%%%%%%%
\section{Strange limits}

We now develop methods that will be useful for identifying the
closures of various Beurling subspace manifolds.

\begin{prop}
  \label{prop:sot-limit}
  Let $B$ be a union of intervals of $\bR$. A sequence of
  projections~$P_n$ converges in the strong operator topology to
  $M_{\chi_B}$ if and only if

  (i) $\|P_n\chi_E\|_2\to\|\chi_E\|_2$
  for every compact interval $E\subseteq B$, and

  (ii) $\|P_n\chi_F\|_2\to 0$
  for every compact interval $F\subseteq\bR\setminus B$.
\end{prop}
\begin{proof}
  Necessity is clear.
  Suppose that (i) and~(ii) hold, and let $E$ and $F$ be disjoint compact intervals. If
  either is a subset of $\bR\setminus B$ then $|\langle P_n
  \chi_E,\chi_F\rangle|\le\|P_n\chi_E\|\,\|P_n\chi_F\|\to 0$ as $n\to
  \infty$. On the other hand, if $E,F\subseteq B$ then for
  $k\in\{0,1,2,3\}$,
  \begin{align*}
  2\re i^k\langle
  P_n\chi_E,\chi_F\rangle&=\|P_n(i^k\chi_E+\chi_F)\|^2-\|P_n\chi_E\|^2-\|P_n\chi_F\|^2
    \\&\le \|\chi_E\|^2-\|P_n\chi_E\|^2 +
    \|\chi_F\|^2-\|P_n\chi_F\|^2
	\to 0
  \end{align*}
  as $n\to \infty$, and we again conclude that $\langle
  P_n\chi_E,\chi_F\rangle\to 0$ as $n\to\infty$.

  Since the unit ball of $\B(\lt)$ is compact and metrisable in the
  weak operator topology,
  there is a positive contraction $C$ which is a weak cluster point of
  $\{ P_n \}$, say $P_{n_k}\to C$ weakly. Now $\langle
  C\chi_E,\chi_F\rangle=0$ for disjoint compact intervals
  $E,F$ which are each contained in either $B$ or its complement; an
  approximation argument shows that this remains the case whenever $E$
  and $F$ are disjoint bounded measurable sets. Approximation by
  simple functions yields $\langle C \chi_E,\chi_F f\rangle=0$ for
  every $f\in\lt$, and similarly $\langle C\chi_E g, \chi_F
  f\rangle=0$ for all $f,g\in\lt$. It follows that $C=M_\phi$ for some
  $\phi\in \Li$ with $0\le\phi\le1$ and $\supp\phi\subseteq
  B$. If $E$ is a compact interval which is contained in $B$, then
  $\|P_{n_k}\chi_E\|^2=\langle P_{n_k}\chi_E,\chi_E\rangle\to \|\chi_E\|^2 =
  \langle\phi \chi_E,\chi_E\rangle$, which forces $\phi$ to be equal
  to $1$ almost everywhere on $E$. So $C=M_{\chi_B}$. Thus the weak
  limit of $\{P_{n_k}\}$ is a projection, $M_{\chi_B}$, from which it
  follows that $\sotlim_{k\to\infty}P_{n_k}=M_{\chi_B}$ as well.

  Thus every subsequence of $\{P_n\}$ has a
  subsubsequence whose limit in the strong operator topology is
  $M_{\chi_B}$. Since the unit ball of $\B(\lt)$ is metrisable in this
  topology, this shows that $P_n\to M_{\chi_B}$ strongly.
\end{proof}

We now determine conditions under which we can identify the strong
operator topology
convergence $P_n\to M_{\chi_B}$ for a sequence of Beurling
projections $P_n=[e^{ik_n\psi}\Ht]$, where $k_n\to\infty$ and $B$
depends on the derivative of $\psi(x)$.

\begin{defn}
  A function $\psi:\bR\to\bR$ is
  \textit{admissible}
  if the following subsets of $\bR$ are discrete (that is, they
  have no  accumulation points):
  \begin{enumerate}
  \item the set $\Gamma(\psi)$ of points at which $\psi$ fails to
    be twice continuously differentiable;
  \item $(\psi')^{-1}(0)$; and
  \item the set $\Lambda(\psi)$ consisting of the points at which
    $\sgn(\psi'')$ is not locally constant.
  \end{enumerate}
\end{defn}

For example, non-constant rational functions and the trigonometric
functions are easily seen to be admissible, as is the map $x\mapsto
\log|x|$.

Let $F:\Lt\to\Lt$ be the unitary Fourier-Plancherel transform
defined by
\[ Ff(x)=\frac1{\sqrt{2\pi}}\int_{-\infty}^\infty f(y)e^{-ixy}\,dy,\quad
f\in \Lt.\]
The Hardy space $\Ht$ is equal to $F^*L^2(0,\infty)$, so
\[[e^{ik\psi}\Ht]=\Ad(M_{e^{ik\psi}}F^*)M_{\chi_{(0,\infty)}}\] where
$\Ad(U)T=UTU^*$ for a unitary operator $U$ and $T\in\B(\H)$, and
$M_\phi$ is the multiplication operator on $\Lt$ corresponding to a
function $\phi\in\Li$, which is unitary if $\phi$ is unimodular.
In particular, if $P=[e^{ik\psi}\Ht]$ then
\begin{align*}
\|P\chi_S\|_2 &=
\|(\Ad(M_{e^{ik\psi}}F^*)M_{\chi_{(0,\infty)}})\chi_S\|_2 \\ &=
\|M_{e^{ik\psi}}F^*(M_{\chi_{(0,\infty)}}FM_{e^{ik\psi}}^*\chi_S)\|_2=
\| \chi_{(0,\infty)}F(e^{-ik\psi}\chi_S)\|_2.
\end{align*}
It is therefore of interest to find an estimate for the Fourier-Plancherel
transform $F(e^{-ik\psi}\chi_S)$ on the positive half-line. This is
done in the next lemma for admissible functions~$\psi$, using an
integration by parts in the spirit
of~\cite{stein93:harmonic-analysis}, Chapter~VIII.

\begin{lem}
  \label{lem:osc}
  Let $\psi$ be admissible and let $S\subseteq\bR$ be a compact
  interval of positive length such that $\Gamma(\psi)\cap S=\emptyset$ and
  $\psi'(S)\subseteq(0,\infty)$. Let $\Delta$ be the function
  $\Delta(z)=2(z+\alpha)^{-1}$ where $\alpha=\min\{\psi'(x):x\in S\}$
  and let $N=|\Lambda(\psi)\cap S^\circ|+2$.  Then for $k>0$ and
  almost every $z>0$, $\sqrt{2\pi}\,|F(e^{-ik\psi}\chi_S)(kz)| \le
  Nk^{-1}\Delta(z)$.
\end{lem}
\begin{proof}
  Since $\Lambda(\psi)$ is discrete,
  finitely many points in $\Lambda(\psi)$ lie in the interior of $S$, say
  $\lambda_2<\lambda_3<\dots<\lambda_{N-1}$. We also write
  $\lambda_1,\lambda_{N}$ for the boundary points of $S$ so that
  $S=[\lambda_1,\lambda_{N}]$. Since $\Gamma(\psi)\cap S=\emptyset$,
  the function $\psi$ is twice continuously differentiable on $S$ and
  for $z>0$,
  \begin{align*}
    \sqrt{2\pi}\,|F&(e^{-ik\psi}\chi_S)(kz)|
    = \Big| \int_S e^{-ik(xz+\psi(x))}\,dx\Big|\\[6pt]
    &= \Big| \int_S \frac{
      i}{k(z+\psi'(x))}\frac{d}{dx}\big(e^{-ik(xz+\psi(x))}\big)\,dx\Big|
    \displaybreak[0]\\[6pt]
    &=\frac1{k}\Big| \Big[  \frac{e^{-ik(xz+\psi(x))}
    }{z+\psi'(x)} \Big]_{\lambda_1}^{\lambda_{N}} -
    \int_S  \frac{d}{dx}\Big(\frac1{z+\psi'(x)}\Big)
    e^{-ik(xz+\psi(x))}\,dx \Big|\\[6pt]
    &\le \frac1{k}
    \Big(
    \Delta(z)+
    \int_S \Big| \frac d{dx}\Big(\frac{1}{z+\psi'(x)}\Big)\Big| \,dx
    \Big).
  \end{align*}
  Let $j\in\{1,2,\dots,N-1\}$. For $x\in(\lambda_j,\lambda_{j+1})$, the quantity
  \[\sgn \frac d{dx}\Big(\frac{1}{z+\psi'(x)}\Big) =
  \sgn\frac{-\psi''(x)}{(z+\psi'(x))^2}=-\sgn\big(\psi''(x)\big)\] is
  constant, say $\sigma_j\in\{-1,0,1\}$.
  So
  \begin{align*}
    \int_S \Big| \frac d{dx}\Big(\frac{1}{z+\psi'(x)}\Big)\Big| \,dx
    &= \sum_{j=1}^{N-1} \sigma_j \int_{\lambda_j}^{\lambda_{j+1}} \frac
    d{dx}\Big(\frac{1}{z+\psi'(x)}\Big) \,dx\\
    &= \sum_{j=1}^{N-1} \sigma_j
    \Big[\frac{1}{z+\psi'(x)}\Big]_{\lambda_j}^{\lambda_{j+1}}\\
    &\le (N-1)\Delta(z).
  \end{align*}
  The result follows.
\end{proof}

\begin{thm}
  \label{thm:psi}
  Let $k_n$ be a sequence of positive numbers with $k_n\to\infty$ as
  $n\to\infty$, let $\psi$ be admissible and let
  $P_n=[e^{ik_n\psi}\Ht]$ for $n\in\bN$. Then $P_n\sotto [L^2(B_-)]$
  as $n\to\infty$, where 
  $B_-=(\psi')^{-1}\big((-\infty,0)\big)$.
\end{thm}
\begin{proof}
  Let $B_+=(\psi')^{-1}\big((0,\infty)\big)$ and
  $B_0=(\psi')^{-1}(0)$.
  Let $S$ be a compact subinterval of $B_+$ of positive length which
  does not intersect $\Gamma(\psi)$.  We will show that
  $P_n\chi_S\to0$.  Since
  $P_n=\Ad(M_{e^{ik_n\psi}}F^*)M_{\chi_{(0,\infty)}}$,
  \begin{align*}
    \|P_n\chi_S\|_2^2 &= \|
    \chi_{(0,\infty)}F(e^{-ik_n\psi}\chi_S)\|_2^2\\
    &= \int_0^\infty |F(e^{-ik_n\psi}\chi_S)(y)|^2\,dy\\
    &= \,k_n\!\int_0^\infty |F(e^{-ik_n\psi}\chi_S)(k_nz)|^2\,dz
  \end{align*}
  where we have made the change of variables $y=k_nz$.
  We apply Lemma~\ref{lem:osc}:
  \[\|P_n\chi_S\|_2^2\le \, \frac{k_n}{2\pi}\int_0^\infty \Big(\frac
  N{k_n}\Big)^2|\Delta(z)|^2\,dz = \frac{N^2}{2\pi k_n}
  \|\chi_{(0,\infty)}\Delta\|_2^2\to 0\]  as $n\to\infty$,
  where $\Delta(z)\in L^2(\chi_{(0,\infty)}\,dz)$ and $N\in\bN$ are
  defined as in Lemma~\ref{lem:osc}.

  Let $\phi=-\psi$, let $Q_n=[e^{ik_n\phi}\Ht]$ and let $C$ be
  the conjugation operator $Cf=\overline f$ for $f\in
  \Lt$. Observe that $Q_n =CP_n^\perp C$. Applying the argument
  of the previous paragraph to $\phi$ instead of $\psi$ shows
  that $Q_n\chi_T\to 0$ whenever $T$ is a compact subinterval of
  $B_-\setminus\Gamma(\psi)$. Hence $P_n\chi_T\to\chi_T$.

  We have shown that $P_n\chi_S\to0$ when $S$ is a compact subinterval
  of $B_+\setminus\Gamma(\psi)=\bR\setminus(B_-\cup B_0\cup
  \Gamma(\psi))$ and $P_n\chi_T\to \chi_T$ when $T$ is a compact
  subinterval of $B_-\setminus\Gamma(\psi)$. Since $\Gamma(\psi)$ and
  $B_0$ are discrete sets, the same conclusion holds for $S$ and $T$
  arbitrary compact subintervals of $\bR\setminus B_-$ and $B_-$
  respectively. Now Proposition~\ref{prop:sot-limit} completes the
  proof.
\end{proof}

The theorem enables us to compute immediately a wide variety of
strange limits $P_n\sotto P$, so-called because while every
nonzero function in the range of $P_n$ has full support, those
for $P$ itself are supported in a proper measurable set.
For example, $[e^{-inx^2/2}\Ht]
\sotto [\Ltp]$ as $n\to \infty$, and $[e^{in(x^3+bx^2+cx)}\Ht]
\sotto \Ltwo{[\alpha,\beta]}$ if the roots $\{\alpha,\beta\}$ of
the equation $3x^2+2bx+c=0$ are real, and has limit~$0$ otherwise.  We also
remark that when $\psi$ is admissible,
\[\sotlim_{n\to\infty}[e^{in\psi}\Ht]=
\big(\!\sotlim_{n\to\infty}[e^{-in\psi}\Ht]\big)^\perp.\]

While admissible functions are adequate for our applications one can
perhaps partially relax this constraint.  However, we are unaware of a
general formula for the limit when $\psi$ is real, measurable and
locally bounded.

The next corollary will play a part in the proof of
Theorem~\ref{thm:sot-closures}.

\begin{cor}
  \label{cor:psi_n}
  Let $k_n$, $\psi$ and $B_\pm$ be as above, let $\psi_n$ be a
  sequence of admissible functions and let $P_n=[e^{ik_n\psi_n}\Ht]$.
  Suppose that the set $\Gamma=\Gamma(\psi)\cup\bigcup_n
  \Gamma(\psi_n)$ is discrete. Let \[ \I=\{
  S\subseteq\bR\setminus\Gamma:S=[\alpha,\beta],\ \alpha<\beta\}\] and
  suppose that the quantity $N(S)=\sup_n|\Lambda(\psi_n)\cap S^\circ|$
  is finite for every $S\in\I$.
    If $\psi_n'\to\psi'$ uniformly on~$S$ for every interval
    $S\in\I$ then $P_n\to [L^2(B_-)]$ strongly as $n\to\infty$.
\end{cor}

\begin{proof} Choose a compact subinterval $S\subseteq
  B_+\setminus\Gamma$.  Let \[\alpha=\min\{\psi'(x):x\in
  S\}\qand\alpha_n=\min\{\psi_n'(x):x\in S\}\text{ for } n\in\bN.\] Pick $n$
  sufficiently large that $\|\psi_n'-\psi'\|_S<\alpha/2$; then
  $\alpha_n>\alpha/2>0$. Writing
  \[\Delta_n(z)=2(z+\alpha_n)^{-1}\qand
  \tilde\Delta(z):=2(z+\alpha/2)^{-1},\]
  Lemma~\ref{lem:osc} applies as before to show that
  \[ \|P_n\chi_S\|^2 \le
  \frac{N(S)^2}{2\pi k_n}\|\chi_{(0,\infty)}\Delta_n\|_2^2\le
  \frac{N(S)^2}{2\pi k_n}\|\chi_{(0,\infty)}\tilde\Delta\|_2^2.\]
  Since $\tilde\Delta(z)\in
  L^2(\chi_{(0,\infty)}\,dz)$, this shows that $P_n\chi_S\to0$.
  The remainder of the proof proceeds as in Theorem~\ref{thm:psi}.
\end{proof}

Beurling's characterisation of invariant subspaces for the
bilateral shift operator when transferred to the setting $\Lt$
amounts to the identification of $\Lat\{M_\lambda :\lambda\ge0\}$
with the disjoint union
\begin{multline*}
  \{u\Ht:\text{$u$ unimodular}\}\cup\{L^2(E):\text{$E$ measurable}\}\\
  =\M_\pure\cup\Lat\{M_\phi:\phi\in\Li\}.
\end{multline*}
Here $\Lat\A$ denotes the lattice of closed invariant
subspaces for a family of operators~$\A$, and $\M_{\pure}$ is the
set of invariant subspaces~$K$ which  are purely invariant in the
sense that the intersection of the subspaces $M_\lambda K$ for
$\lambda\ge0$ is
trivial. We now use  the methods of this section to
show that
$\Lat\big(\Li\big)\subseteq\clo{\M_\pure}$. This seems to be a previously
unobserved feature in the classical setting which may well have a wider
manifestation. However,
the authors are unaware of any general results of this nature.

\begin{lem}
  \label{lem:noacc}
  Let $m$ denote Lebesgue measure on~$\bR$.
  If $B\subseteq\bR$ is measurable and $\epsilon>0$ then there is a
  countable disjoint union of open intervals $V$ such that $\bd V$ is
  discrete and $m(V\sd B)<\epsilon$.
\end{lem}
\begin{proof}
  Fix $n\in\bZ$ and write $B_n=B\cap(n,n+1)$ and $\epsilon_n=
  2^{-|n|}\epsilon/3$.
  Using elementary properties of Lebesgue measure,
  we can find a set $U_n=\bigcup_{j\ge1} I_j\supseteq B_n$ such
  that $\{I_j\}_{j\ge1}$ are disjoint open subintervals of
  $(n,n+1)$ and $m(U_n\setminus B_n)<\frac12\epsilon_n$. Pick $k$ so
  that $\sum_{j>k} m(I_j) < \frac12\epsilon_n$ and let
  $V_n=\bigcup_{1\le j\le k}I_j$. Now
  \[ V_n\sd B_n = \big( (U_n\setminus B_n)\cap V_n\big)\cup
  (B_n\setminus V_n) \subseteq U_n\setminus B_n \cup\bigcup_{j>k} I_j,\]
  so $m(V_n\sd B_n)< \tfrac12 \epsilon_n + \tfrac12\epsilon_n =
  \epsilon_n$.

  Repeat for each $n\in\bZ$ and let $V=\bigcup_{n\in\bZ}V_n$.
  Then \[m(V\sd B)=\sum_{n\in\bZ}m(V_n\sd B_n)<
  \sum_{n\in\bZ}\epsilon_n=\epsilon\] and $\bd V$ is discrete, since
  $\bd V\cap [n,n+1]$ is finite for each $n$.
\end{proof}

\begin{thm}\label{thm:purelimits}
  If $B$ is any measurable subset of $\bR$ then there is a sequence of
  projections $P_n=[e^{ik_n\psi_n}\Ht]$ where
  $k_n>0$ and each $\psi_n$ is a real-valued function such that
  $\sotlim_{n\to\infty}P_n=M_{\chi_B}$.
\end{thm}
\begin{proof}
  Let $\epsilon>0$. By Lemma~\ref{lem:noacc}, we can find a countable
  disjoint union of open intervals $V_\epsilon$ such that $\bd
  V_\epsilon$ is discrete and $m(V_\epsilon\sd
  B)<\epsilon$. The function $\psi_\epsilon(x)= x(\chi_{\bR\setminus
    V_\epsilon} - \chi_{V_\epsilon})$ satisfies
  $\Gamma(\psi_\epsilon)\cup\Lambda(\psi_\epsilon)\subseteq\bd
  V_\epsilon$ and $(\psi_\epsilon')^{-1}(0)=\emptyset$, so
  $\psi_\epsilon$ is admissible. Let
  $P_{\epsilon,n}=[e^{in\psi_\epsilon}\Ht]$. By
  Theorem~\ref{thm:psi},
  $\sotlim_{n\to\infty}P_{\epsilon,n}=M_{\chi_{V_\epsilon}}$, and we
  also have $M_{\chi_{V_\epsilon}}\to M_{\chi_B}$ strongly as $\epsilon \to 0$.

  Let $d$ be a metric inducing the strong operator topology on the
  unit ball of $\B(\lt)$ and let $n\in\bN$. Choose $\epsilon_n>0$ such
  that $d(M_{\chi_{V_{\epsilon_n}}}, M_{\chi_B})<1/2n$ and then choose
  $k\in\bN$ so that $P_n:=P_{\epsilon_n,k}$ satisfies
  $d(P_n,M_{\chi_{V_{\epsilon_n}}})<1/2n$. Now
  $d(P_n,M_{\chi_B})<1/n$ so $P_n\to M_{\chi_B}$ strongly as
  $n\to\infty$.
\end{proof}

%%%%%%%%%%%%%%%%%%%%%%%%%%%%%%%%%%%%%
\section{Closures of Beurling subspace manifolds}
\label{sec:4}

We now obtain sufficient conditions under which Beurling subspace
manifolds $\M(\S)$ have closures, in the strong operator topology,
which are compact. Using this we construct various $n$-spheres and
$n$-balls in $\Proj(\Lt)$. At the end of the section we pose
some further lines of enquiry.
\medskip

Let $f=(f_1,f_2,\dots,f_n)$ be an $n$-tuple of functions
$f_j:\bR\to\bR$. We write $\langle
f,\lambda\rangle=\lambda_1 f_1+\lambda_2
f_2+\dots + \lambda_n f_n$ for $\lambda\in\bR^n$, and \[\S_f=\{\langle
f,\lambda\rangle:\lambda\in\bR^n\}.\]
\begin{defn}
  The $n$-tuple
  $f$ is \emph{admissible} if
  \begin{enumerate}
  \item the set $\{f_1,f_2,\dots,f_n\}$ is linearly independent over~$\bR$;
  \item every nonzero function in $\S_f$ is admissible; and
  \item $\sup_{g\in\S_f\setminus\{0\}}|K\cap\Lambda(g)|<\infty$ for each compact
    set $K\subseteq\bR$.
  \end{enumerate}
  We will also write $\Gamma(f)=\bigcup_{j=1}^n \Gamma(f_j)$ and
  remark that this is equal to $\bigcup_{g\in
    \S_f}\Gamma(g)$ and is plainly discrete.
\end{defn}

Given an admissible $n$-tuple $f$, let $\theta:\bR^n\to\Proj(\Lt)$
be the map $\lambda\mapsto [e^{i \langle f,\lambda\rangle}\Ht]$.
Observe that $\theta$ is strongly continuous, since if
$\lambda\ub k\to \lambda$ in $\bR^n$ then $\langle f, \lambda\ub
k\rangle \to \langle f,\lambda\rangle$ uniformly on compact
subsets of $\bR\setminus \Gamma(f)$ and so \[\theta(\lambda\ub k)
= \Ad(M_{e^{i\langle f,\lambda_k\rangle}})[\Ht]\sotto
\Ad(M_{e^{i\langle f,\lambda\rangle}})[\Ht]=\theta(\lambda)\]
as $k\to\infty$.

We write $\M(\S_f)=\theta(\bR^n)$. We will shortly see that $\M(\S_f)$
is a Beurling subspace manifold.

\begin{thm}
  \label{thm:sot-closures}
  Given an admissible $n$-tuple~$f$, the
  closure of the range of $\theta$ in the strong operator topology is
  \[ \overline{\M(\S_f)} = \theta(\bR^n) \cup
  \big\{ [L^2\big((\psi')^{-1}((-\infty,0)\big)] :
    \psi\in\S_f\setminus\{0\}\big\}.
  \]
\end{thm}
\begin{proof}
  Let $\lambda\ub k$ be a sequence in $\bR^n$ and let
  $P_k=\theta(\lambda^{(k)})$ be the corresponding sequence of
  projections. Passing to a subsequence, we may assume that
  $\lambda\ub k$ converges to a vector $\lambda\in
  (\bR\cup\{\pm\infty\})^n$ as $k\to\infty$.  If $\lambda$ actually
  lies in $\bR^n$ then $P_k\to \theta(\lambda)$ by
  the continuity of $\theta$.
  Otherwise, if $\alpha_k=\max_j|\lambda_j\ub k| $ then $\alpha_k\to
  \infty$ as $k\to\infty$. Let $\mu\ub k=\alpha_k^{-1}\lambda\ub k$.
  Passing to a subsequence, we may assume that $\alpha_k =
  |\lambda_{j_0}\ub k|$ for some $j_0$ independent of $k$, and that
  $\mu\ub k\to \mu$ for some $\mu\in[-1,1]^n$. Since $\mu_{j_0}=\pm1$,
  this limit $\mu$ is nonzero. Now $P_k=[e^{i\alpha_k \psi_k}\Ht]$ where
  $\psi_k=\langle f_k, \mu\ub k\rangle$, and if $\psi=\langle
  f,\mu\rangle$ then $\psi_k'\to \psi'$ uniformly on compact subsets
  of $\bR\setminus \Gamma(f)$. By
  Corollary~\ref{cor:psi_n}, $P_k\to
  [L^2((\psi')^{-1}((-\infty,0)))]$.

  Conversely, if $\psi\in\S_f\setminus\{0\}$ then by
  Theorem~\ref{thm:psi},
  \[
  [L^2\big((\psi')^{-1}((-\infty,0))\big)]=\sotlim_{n\to\infty}[e^{in\psi}\Ht]\in\overline{\M(\S_f)}.
  \qedhere
  \]
\end{proof}

\begin{defn}
  \label{defn:sim}
  Given an admissible $n$-tuple $f$, let $\sim$ be the equivalence
  relation defined on $\bnc$ by $\lambda\sim\mu$ if $\lambda=\mu$ or
  \[ \lambda,\mu\in\sn\qand m(\{x: \langle f',\lambda\rangle(x)>0\} \sd
  \{x: \langle f',\mu\rangle(x)>0\})=0. \]
  Here $f'=(f_1',\dots,f_n')$ and $m$ is Lebesgue measure on~$\bR$.
  We write $\bmodsim$ for
  the corresponding topological quotient space.
\end{defn}

\begin{prop}
  \label{prop:compact-quotients}
  The topological space $\bmodsim$ is homeomorphic to
  $\clo{\M} = \clo{\M(\S_f)}$. In particular, $\clo\M$ is compact.
\end{prop}
\begin{proof}
  Let $\theta$ be the continuous map $\bR^n\to\clo\M$ defined
  above. Observe that $\theta$ is injective: for
  $\theta(\lambda)=\theta(\mu)$ if and only if
  $e^{i\langle f,\lambda-\mu\rangle} \Ht=\Ht$ which implies that the
  function $\langle f, \lambda-\mu\rangle$ is constant modulo
  $2\pi$ almost everywhere. Such a function cannot be admissible,
  so it is zero; by linear independence, $\lambda=\mu$.

  Let $\alpha:\bR^n\to\bB^n$ be a homeomorphism of the form
  \[\alpha:\lambda\mapsto \rho(\|\lambda\|) \lambda\|\lambda\|^{-1}\] where
  $\rho:[0,\infty)\to[0,1)$ is a homeomorphism. Consider the injective
  continuous map $\phi=\theta\circ \alpha^{-1}:\bB^n\to \clo\M$. We
  extend this to $\bnc$ by defining $\phi(\lambda)=\lim_{r\uparrow
    1}\phi(r\lambda)$ for $\lambda\in\sn$; this limit exists by
  Theorem~\ref{thm:psi}. The extended map is also continuous by
  Corollary~\ref{cor:psi_n} and surjective by
  Theorem~\ref{thm:sot-closures}. Since $\phi(\lambda)=\phi(\mu)$ if
  and only if $\lambda\sim\mu$, it follows that $\phi$ induces a homeomorphism from
  the compact space $\bmodsim$ to the Hausdorff space $\clo\M$.
\end{proof}

\begin{rem}
  This proof shows that $\theta:\bR^n\to\M(\S_f)$ is a
  homeomorphism when $f$ is admissible, and so $\M(\S_f)$ is indeed a
  Beurling subspace manifold.
\end{rem}

Determining the precise nature of the quotient space $\bmodsim$ seems
difficult in general. However, the following shows that there are no
surprises when $n=2$.

\begin{lem}
  Let $\I$ be any set of closed, pairwise disjoint intervals, each
  contained in $(0,1)$. Then there exists a continuous increasing
  surjection $\beta:[0,1]\to[0,1]$ such that $\beta(x)=\beta(y)$ if
  and only if $x,y\in I$ for some $I\in\I$.
\end{lem}
\begin{proof}
  We may assume without loss of generality that each $I\in \I$ is of
  the form $I=[a,b]$ with $0<a<b<1$. Observe that $\I$ must be
  countable since for each $n\ge1$, the set
  $\{[a,b]\in\I:b-a>n^{-1}\}$ is finite.  Enumerate $\I$ in decreasing
  size, so that $\I=\{I_1,I_2,\dots\}$ where $|I_n|\ge|I_{n+1}|$ for
  $n\ge1$. Assign values $a_n$ on $I_n$ to create a partially defined
  increasing function $\beta_1$ on $\bigcup\I$ as follows.  Set
  $\beta_1(0)=0$ and $\beta_1(1)=1$. Then set $a_1=\tfrac12$ and
  $a_2=\tfrac14$ or $\tfrac34$, according to whether the new interval
  $a_2$ is to the left or the right of $a_1$. Continue
  ``interpolating dyadically'' in this way, so that $a_{n+1}$ is
  chosen as the mean of $\beta_1(\ell)$ and $\beta_1(r)$ where $\ell$
  (respectively,~$r$) is the point of $\{0,1\}\cup I_1\cup\dots\cup I_n$
  immediately to the left (respectively, right) of $I_{n+1}$.

  Let $K$ be the closure in $[0,1]$ of $\{0,1\}\cup\bigcup\I$. Let
  $U$ be the complement of $K$, an open set in $[0,1]$ and thus a
  union of open intervals, which we call ``gaps''. We call a gap
  ``good'' if it is of the form $J=(b,a')$ where $I_n=[a,b]$ and
  $I_m=[a',b']$ are intervals in $\I$, and ``bad'' otherwise. We can
  extend $\beta_1$ to a good gap $J$ by linear interpolation between
  $a_n$ and $a_m$.  On the other hand, a bad gap $J'$ must have at
  least one of its end points a limit of endpoints of intervals
  $I\in\I$, so $\beta_1$ can be extended to the closure $\overline{J'}$ as a
  non-decreasing function in a unique way, namely by being constant on
  $\overline{J'}$. So we get a continuous increasing function $\beta_1$ defined
  on $[0,1]$.  We now correct the constancy on the intervals $\overline{J'}$ by
  forming $\beta_3 = \beta_1 + \beta_2$ where
  \[\beta_2 = \sum_{\text{bad $J'$}} \beta_{J'}\]
  and $\beta_{J'}$ is the continuous map on $[0,1]$ which takes the
  value $0$ to the left of $J'$, the value $|J'|$ to the right and is
  linear on $J'$. (None of these $\beta_{J'}$ spoil constancy on the
  $I_n$.) Continuity of $\beta_3$ is clear. Finally, let $\beta =
  c\beta_3$ with $c=(1+\sum_{\text{bad $J'$}}|J'|)^{-1}$ so that
  $\beta(1) = 1$.
\end{proof}

\begin{lem}
  \label{lem:circle-quotients}
  Let $\I$ be a set of closed, pairwise disjoint proper arcs
  in $\bT$. Then $\cd/\I$ is homeomorphic to~$\cd$.
\end{lem}
\begin{proof}
  We may assume that $(1,0)$ does not belong to any set in $\I$.
  Transfer $\I$ to $[0,1]$ in the obvious manner and
  apply the previous lemma to obtain an increasing continuous
  surjection $\beta:[0,1]\to [0,1]$ such that $\beta(x)=\beta(y)$ if
  and only if $\{e^{2\pi i x}, e^{2\pi i y}\}\subseteq I$ for some
  $I\in\I$.  The map $\phi:\cd\to\cd$ given by
  $\phi(re^{2\pi ix})=re^{2\pi i(r\beta(x)+(1-r)x)}$ for $r,x\in (0,1]$
  and $\phi(0)=0$ induces a continuous bijection $\cd/\I\to\cd$, which
  is a homeomorphism since $\cd/\I$ is compact.
\end{proof}

\begin{prop}
  \label{prop:admissible-pairs}
  For every admissible pair $(f,g)$, the closure in the strong
  operator topology of the Beurling subspace manifold
  $\M(\S_{(f,g)})$ is homeomorphic to~$\cd$.
\end{prop}
\begin{proof}
  Let $\sim$ be the equivalence relation on $\cd$ of
  Definition~\ref{defn:sim}, with equivalence classes
  $\{[\lambda]:\lambda\in\cd\}$.  Observe that
  $\lambda\not\sim-\lambda$ and $\lambda\sim\mu$ if and only if
  $-\lambda\sim-\mu$. Also, the equivalence classes are connected: if
  $\lambda\sim\mu$ then the shorter of the two arcs joining $\lambda$
  to $\mu$ is contained in $[\lambda]=[\mu]$.

  We show that the equivalence classes are also closed.  If there is a
  non-trivial equivalence class, then we can make a different choice
  of $f$ and $g$ without changing the set $\S_{(f,g)}$ to arrange that
  $(1,0)\sim(0,-1)$, and so also that $(-1,0)\sim(0,1)$. Any remaining
  equivalence classes are of the form $[\lambda]$ where
  $\lambda=(\lambda_1,\lambda_2)$ with $\lambda_1\lambda_2>0$; we may
  assume that $\lambda_1>0$ and $\lambda_2>0$. Let $h=g'/f'$ and let
  $\alpha=-\lambda_1/\lambda_2$. Then
  $[\lambda]=[\lambda]_+\cap[\lambda]_-$ where
  \begin{align*}
    [\lambda]_+ &= \big\{ \mu :
    \{f'>0\}\cap(\{ h>\alpha \} \sd \{ h>-\mu_1/\mu_2\} )
    \text{ is null} \}\text{ and} \\
    [\lambda]_- &= \big\{ \mu :
    \{f'<0\}\cap(\{ h<\alpha \} \sd \{ h<-\mu_1/\mu_2\} )
    \text{ is null} \}.
  \end{align*}
  Here and below we employ abbreviations of the form
  $\{P(\phi)\}$ to mean the set $\{x\in\bR:P(\phi(x))\}$ where $\phi:\bR\to\bR$ and $P$
  is a predicate depending on a real parameter.

  Note that for any constant $k$, the set $\{x:h(x)=k\}$ is null; for
  if not, then $f'-kg'$ takes the value $0$ on a non-null set, so
  $f-kg$ is not admissible, whence $f=kg$; but $f$ and $g$ are
  independent.
  For $\gamma>0$, let $u_\gamma$ be the unit vector
  $(\gamma,1)/(\gamma^2+1)^{1/2}$.  Each $[\lambda]_+$ is clearly
  connected. It is also closed, for if $b>a>0$ and
  $\{u_\gamma:\gamma\in(a,b)\}\subseteq[\lambda]_+$ then
  $\{\alpha<h\le -b\}\cap\{f'>0\}$ is equal to the union of the
  null sets $\{f'>0,\ h=-b\}$ and $\bigcup_{n\ge1}\{\alpha < h \le
  -b-n^{-1}\}\cap\{f'>0\}$, so is null; hence $u_b$, and by a similar
  argument $u_a$, lie in $[\lambda]_+$. The class $[\lambda]_-$ is
  closed for the same reasons, and hence each class $[\lambda]$ is a
  closed subarc of $\bT$.  Now
  Proposition~\ref{prop:compact-quotients} and
  Lemma~\ref{lem:circle-quotients} complete the proof.
\end{proof}

Given an admissible $n$-tuple $f$, let $\Sigma(\S_f)$ denote the
\emph{complement completion} of $\M(\S_f)$; that is, the
closure of $\M(\S_f)\cup\M(\S_f)^\perp$ in the \sot. It is not
hard to see that $\M(\S_f)$ and $\M(\S_f)^\perp$ are disjoint, and the
boundaries of these sets are equal by
Theorem~\ref{thm:sot-closures}. From this it follows that
provided $\overline{\M(\S_f)}$ is homeomorphic to
$\overline{\bB^n}$, the set $\Sigma(\S_f)$ is homeomorphic
to $S^n$ since it is homeomorphic to the union of two copies of
$\overline{\bB^n}$ joined at their boundaries.

\begin{eg}
  \label{eg:FPS}
  The \emph{Fourier-Plancherel sphere} is the set of projections
  \[\FPS=\Sigma(\S_f) \quad \text{where $f$ is the admissible pair $f=(x,x^2)$}.\] By
  Proposition~\ref{prop:admissible-pairs} or the direct arguments
  of~\cite{pow-kat:FB}, we see that $\overline{\M(\S_f)}$ is
  homeomorphic to $\overline\bD$, and as observed
  in~\cite{pow-kat:hyp}, the order structure of $\FPS$ is that of a
  union of continuous nests which meet only at $0$ and
  $I$ and $\FPS$ is homeomorphic to a 2-sphere on
  which the Fourier-Plancherel transform $F$ acts as a
  quarter-rotation. In particular, $\FPS\setminus\{0,I\}$ is a locally
  unitary subspace manifold: we clearly have a locally unitary
  structure on $\M(\S_f)\cup \M(\S_f)^\perp$, and we can use $F$ to
  transfer this structure to the remaining subspaces.
\end{eg}

\begin{eg}
  \label{eg:2dim}
  Let $f$ be the admissible pair $f=(x^{-1},x)$. An easy extension
  of~\cite{pow-kat:hyp}, Lemma~5.1 shows that every projection
  $P\in\M(\S_f)$ lies in a continuum of non-commuting continuous nests
  which intersect only in $\{0,P,I\}$. A simple calculation reveals that the
  equivalence relation $\sim$ has only two non-trivial equivalence
  classes and these are antipodal closed quarter-circles,
  so $\bmodsim[\overline{\bB^2}]$ is homeomorphic to
  $\overline{\bB^2}$ and the boundary projections are
  \[\{P,P^\perp:P=[L^2((-a,a))],\ a\in[0,\infty]\}.\]
  Now although $\Sigma(\S_f)$ is homeomorphic
  to a 2-sphere, $\Sigma(\S_f)\setminus\{0,I\}$ is not
  locally unitary. For if there were a strong operator topology
  neighbourhood of $P_0=[L^2((-a,a))]$ of the form
  \[\N=\{[\rho(\lambda)P_0\Lt]:\lambda\in\bB^2\}\] for some
  unitary-valued representation $\rho$ of $\bR^2$, then $\N$ would
  intersect $\M(\S_f)$ and so contain a projection
  $P=[\rho(\lambda)P_0\Lt]$ for some point $\lambda\in\bB^2$ such that every
  neighbourhood of $P$ contains two non-commuting projections which
  are comparable with $P$.  Applying $\rho(-\lambda)$, we see that
  $\N$ must contain two non-commuting projections which are comparable
  with $P$. However, all the projections in $\Sigma(\S_f)$ which are
  comparable with $P$ all commute since they are of the form $[L^2(E)]$
  for some $E\subseteq\bR$.
\end{eg}

\begin{eg}
  If we take $f=(\tfrac13 x^3,\tfrac12 {x^2}, x)$ then $f'=(x^2,x,1)$
  and it is not hard to check that the corresponding equivalence
  relation $\sim$ on $\overline{\bB^3}$ has two non-trivial
  equivalence classes:
  \begin{align*}
    \{(0,0,1)\}&\cup\{(a,b,c)\in S^2:b^2\le 4ac,\ a>0\} \qand \\
    \{(0,0,-1)\}&\cup\{(a,b,c)\in S^2:b^2\le 4ac,\ a<0\}
  \end{align*}
  which correspond to $0$ and $I$ respectively when we identify the
  quotient space $\bmodsim[\overline{\bB^3}]$ with
  $\overline{\M(\S_f)}$. These equivalence classes are closed and so
  $\overline{\M(\S_f)}$ is homeomorphic to $\overline{\bB^3}$.  If $P$
  is a projection in $\M(\S_f)$ then we
  claim that the projections $[e^{ik x}P\Lt]$ for $k\in\bR$ are the
  only non-trivial elements of $\overline{\M(\S_f)}$ which are
  comparable with $P$. To see this, recall that no proper subspace of
  the form $[L^2(E)]$ is comparable with $P$ by the F.~\&~M.~Riesz
  theorem, and if $e^{i\langle f,\lambda\rangle}\Ht\subseteq
  e^{i\langle f,\mu\rangle}\Ht$ and $\lambda\ne\mu$ then $e^{i\langle
    f,\lambda-\mu\rangle}\Ht\subseteq\Ht$ and so $e^{i\langle
  f,\lambda-\mu\rangle}$ is a nonzero continuous inner function. This
  must be of the form $\alpha e^{i\beta x}$ for
  a unimodular constant $\alpha$ and $\beta\in\bR$
  (see~\cite{garnett07:bdd-anal-fns}), which verifies the
  claim. On the other hand, the boundary of $\overline{\M(\S_f)}$
  consists of the projections $[L^2(E)]$ where $E$ is either an
  interval or the complement of an interval. As in the previous
  example, it follows that the topological 3-sphere $\Sigma(\S_f)$
  cannot be locally unitary away from $\{0,I\}$ since the local order
  structure changes on the boundary of $\overline{\M(\S_f)}$.
\end{eg}

\begin{eg}
  \label{eg:hyp-sphere}
  Let $f=(x,\log|x|,-x^{-1})$. Then $f'=(1,x^{-1},x^{-2})$ so if
  $\lambda=(a,b,c)$ then $\sgn \langle
  f',\lambda\rangle=\sgn(ax^2+bx+c)$. The equivalence relation
  for $f$ on $\overline{\bB^3}$ is therefore identical to the
  relation considered in the previous example, and
  $\overline{\M(\S_f)}$ is again homeomorphic to
  $\overline{\bB^3}$. The order structure differs however, since
  $\overline{\M(\S_f)}$ contains the set
  $\overline{\M(\S_{(x,x^{-1})})}$ from Example~\ref{eg:2dim}. We
  call the set $\Shyp=\Sigma(\S_f)$ the \emph{hyperbolic
    sphere}. This was first considered
  in~\cite{pow-kat:hyp}, Section~7. We remark that as in
  Example~\ref{eg:2dim}, $\Shyp\setminus\{0,I\}$ cannot be locally unitary.

  We can now easily establish the compactness of the ``extended hyperbolic
  lattice'' $\Mhyp$ considered in~\cite{pow-kat:hyp}; this fact
  was alluded to but not proven there.  For $\alpha\in\bT$, let
  $u_\alpha:\bR\to\bC$ be the two-valued function taking the
  value $1$ on $[0,\infty)$ and $\alpha$ on $(-\infty,0)$. Then
  $\Mhyp$ may be succinctly described as the set of projections
  \[ \Mhyp=\bigcup_{\alpha\in\bT} \Ad(M_{u_\alpha})\Shyp.\] Now
  $\Shyp$ is homeomorphic to the compact space $S^3$ and $\Mhyp$
  is a continuous image of $\bT\times\Shyp$ which is
  compact, so $\Mhyp$ is also compact.
\end{eg}

\begin{eg}
  The equivalence relation $\sim$ need not have a finite number
  of nontrivial equivalence classes. For example, consider the
  admissible triple
  $f=(\tfrac12 x^2,\log|x|,-x^{-1})$ so that if $\lambda=(a,b,c)$ then
  $\sgn\langle f',\lambda\rangle=\sgn(ax^3+bx+c)$. A
  simple analysis of cases reveals that the nontrivial
  equivalence classes for $\sim$ in the upper hemisphere $\{(a,b,c)\in
  S^2:a\ge0\}$ are
  \begin{equation*}
    I_t=\{\lambda/\|\lambda\|: \lambda=(1,s,-t(s+t^2)),\ s\ge
    -3t^2/4\}\cup\{\mu_t\},\quad t\in\bR
  \end{equation*}
  where $\mu_t=(1+t^2)^{-1/2}(0,1,-t)$. The
  class $I_t$ corresponds to the projection
  $[L^2(-\infty,t)]\in\overline{\M(\S_f)}$. Since points of the
  form $(1,s,-t(s+t^2))$ for $s\in\bR$ form a straight line in
  the plane $\{(1,b,c):b,c\in\bR\}$ it follows that $I_t$ is the
  geodesic on $S^2$ joining $\lambda_t/\|\lambda_t\|$ to $\mu_t$
  where $\lambda_t=(1,-3t^2/4,-t^3/4)$. In the $a\le0$ hemisphere
  the nontrivial equivalence classes are the sets $-I_t$ which
  correspond to $[L^2(t,\infty)]$.  It is easy to see that the
  quotient space $\bmodsim[\overline{\bB^3}]$ is homeomorphic to
  $\overline{\bB^3}$; indeed, we may choose a homeomorphism which
  contracts each geodesic $I_t$ to the point $\mu_t$ and extend
  this to all of $\overline{\bB^3}$ is a straightforward manner.
  % a better way to say that succinctly would be nice
  We again conclude that $\overline{\M(\S_f)}$ is homeomorphic to
  $\overline{\bB^3}$.
\end{eg}

\begin{rem}
  We do not know if $\overline{\M(\S_f)}$, or equivalently $\bmodsim$,
  is homeomorphic to $\overline{\bB^n}$ for every admissible $n$-tuple
  $f$. It is natural to try to emulate the argument of
  Proposition~\ref{prop:admissible-pairs}, and it is not hard to show
  that the equivalence classes $[\lambda]$ of $\sim$ satisfy the following
  conditions:
  \begin{enumerate}
  \item $[\lambda]=\{\lambda\}$ if $\lambda\in\bB^n$;
  \item  $[-\lambda]=-[\lambda]$ and
    $[-\lambda]\cap[\lambda]=\emptyset$ for $\lambda\in S^{n-1}$;
  \item if $\lambda\sim \mu$ then $\lambda\sim \nu$ for every $\nu$ on the geodesic in
    $S^{n-1}$ joining $\lambda$ to $\mu$; and
  \item if $\lambda_n\sim \mu_n$ where $\lambda_n\to \lambda$ and $\mu_n\to \mu$ are
    convergent sequences in $S^{n-1}$, then $\lambda\sim \mu$.
  \end{enumerate}
  However, for $n>2$ we have been unable to identify the quotient
  space $\bmodsim$ for such an equivalence relation.
\end{rem}

\begin{rem}
  The order structure of the 2-sphere $\FPS$ and the 3-sphere $\Shyp$
  can be viewed as providing an inherent foliation. We exploit
  this structure in the next section in the determination of their
  unitary automorphism groups. On the other hand we see  in
  Lemma 5.8 that the 2-spheres
  determined by monomial pairs $x^p, x^q$, for $|p|,|q| > 1$ have a
  trivial order structure supported in the common boundary
  of $\M$ and $\M^\bot$.

  The Fourier-Plancherel sphere seems to be a particularly
  distinguished example amongst these 2-spheres.
  Furthermore its equator, $\Sigmae$
  yields an interesting  compact $1$-dimensional subspace manifold
  which is locally unitary and is probably not (periodically) unitary.
  It would be interesting to determine
  other (unitarily inequivalent) subspace manifolds of this form.
\end{rem}

\section{Unitary automorphisms and isomorphisms}
\label{sec:5}

Given a set of projections $\P\subseteq\Proj(\H)$, the unitary
automorphism group of $\P$ is
\[\U(\P)=\{U\in\Unit(\H):(\Ad U)\P = \P\}.\]
As usual, $\Ad U:\B(\H)\to\B(\H)$ is the map $(\Ad
U)T=UTU^*$.
%With a view to a possible classification scheme for
%subspace manifolds and their complement completions, it is of
%interest to compute these unitary invariants.
In this section we compute the unitary automorphism groups of the
Fourier-Plancherel sphere, the hyperbolic sphere and the extended
hyperbolic lattice. These
%complement completions
projection manifolds
inherit a relatively rich order structure from
$\Proj(\Lt)$ which we are able to exploit. In contrast we show in
Section~\ref{sec:5.3} that many other polynomial 2-spheres are essentially
rigid.
%The discussion below is self-contained.
Further operator
 algebra related to the two main examples can be found
in~\cite{power98:repns}, \cite{power06:invariant-subspaces}, \cite{levene-power03:ah-refl}.

\subsection{The Fourier-Plancherel sphere}
\label{sec:5.1}

Recall the definition of the Fourier-Plancherel sphere~$\FPS$ from
Example~\ref{eg:FPS}. The following notation from~\cite{pow-kat:FB} is
convenient:
\begin{gather*}
  \phi_s\in\Li,\quad \phi_s(x)=e^{-isx^2/2},\quad s\in \bR,\\
  V_t\in\Unit(\Lt),\quad V_tf(x)=e^{t/2}f(e^tx),\quad t\in\bR,\\
  M_\lambda\in\Unit(\Lt),\quad M_\lambda=M_{e^{i\lambda x}},\quad
  \lambda\in\bR.
\end{gather*}
We also introduce notation for some nests contained in
$\FPS$:
the \emph{analytic nest} is
$\analnest=\{[M_\lambda\Ht]:\lambda\in\bR\}\cup\{0,I\}$ and the
\emph{Volterra nest} is
\[\volterranest=\{[L^2(t,\infty)]:t\in\bR\}\cup\{0,I\}.\] If we
write
$\N_s=\Ad(\Mphis)\anest$ then
\[
\FPS=\vnest\cup\vnest^\perp\cup\bigcup_{s\in\bR}(\N_s\cup\N_s^\perp),
\]
and the order structure that $\FPS$ inherits from $\Proj(\Lt)$ is such
that if $P,Q\in\FPS\setminus\{0,I\}$ with $P\ne Q$ then $P\vee Q=I$
and $P\wedge Q=0$ unless $\{P,Q\}\subseteq \N$ for some nest $\N$ in
this union.

It is easy to see that $\Mphis$, $M_\lambda$ and $V_t$
all lie in $\U(\FPS)$ for $s,t,\lambda\in\bR$, as does the
Fourier-Plancherel transform $F$ since
\begin{equation}
  \label{eq:fourier}
  (\Ad F)\anest=\vnest,\quad (\Ad F)\N_s=
  \begin{cases}
    \N_{-1/s}^\perp\quad&s>0,\\
    \N_{-1/s}\quad&s<0\\
  \end{cases}
\end{equation}
by~\cite{pow-kat:hyp}, Theorem~7.1. We first show that $V_t$ may
be expressed solely in terms of $\{\Mphis:s\in\bR\}$ and $F$.

\begin{lem}
  \label{lem:Lang-Vt}
  For $t\in \bR$, the dilation operator $V_t$ lies in the
  group generated by $\{\Mphis, F, e^{i\psi}I: s,\psi\in \bR\}$.
  In fact, \[V_t=e^{i\pi/4} \Mphis[ \exp(t)]F  \Mphis[ \exp(-t)]F  \Mphis[ \exp(t)]F.\]
\end{lem}
\begin{proof}
  Let us write $S_{g}$ for the operation of convolution with a
  function $g\in \Li$, defined on the Schwartz~space $\S(\bR)$; that
  is,
  \[ S_{g} f(x)=\int_\bR g(x-t)f(t)\,dt,\quad f\in \S(\bR),\ x\in\bR.\]
  For $\zeta\in\bC\setminus\bR_-$, let $\zeta^{\pm1/2}$ denote
  the square root of~$\zeta^{\pm1}$ with non-negative real part.
  Let $\tilde F$ be the alternate Fourier transform defined on
  $\S(\bR)$ by
  \[ \tilde Ff(x)=\int_\bR f(y)e^{-2\pi ixy}\,dy,\quad f\in\S(\bR),\ x\in\bR.\]
  Observe that $\tilde F=V_{\log 2\pi}F|_{\S(\bR)}$ and that $V_t
  \Mphis = \Mphis[e^{2t}s]V_t$.

  In Section~XI.1
  of~\cite{lang85:sl2(R)} it is shown that
  \begin{equation*}
    \tilde FS_{\phis[2\pi b]}=(ib)^{-1/2}M_{\phis[-2\pi/b]}\tilde F,\quad
b\in\bR\setminus\{0\},
  \end{equation*}
  or, writing $s=2\pi b$ and rearranging,
\begin{equation*}
    S_{\phis}f=(2\pi/is)^{1/2} F^* \Mphis[-1/s]Ff,\quad
    f\in \S(\bR),\ s\in \bR\setminus\{0\}.
  \end{equation*}
  Observe that $\phis(x-t)=e^{isxt}\phis(x)\phis(t)$ for $x,s,t\in\bR$.
  Hence for $x\in \bR$, $s<0$ and $f\in\S(\bR)$,
  \begin{align*}
    S_{\phis}f(x) &= \int_\bR \phis(x-t)f(t)\,dt\notag\\
    &=\phis(x)\int_\bR e^{isxt} \phis(t) f(t)\,dt\notag\\
    &= (2\pi)^{1/2}\phis(x) (F^* \Mphis f)(sx)\\
    &= \big(2\pi/(-s)\big)^{1/2} \Mphis V_{\log (-s)}F\Mphis f(x).
  \end{align*}
  Equating these expressions for $S_{\phis}f$ and using the density of
  $\S(\bR)$ in $\Lt$ gives
  \[
  V_{\log(-s)}=e^{i\pi/4}\Mphis[-s]F^*\Mphis[-1/s]F\Mphis[-s]F^*.
  \]
  Now $F^*=F^3$ and $F^2$ commutes with $\Mphis[\sigma]$ for any
  $\sigma\in\bR$, since $\phis[\sigma]$ is even and $F^2f(x)=f(-x)$.
  So \[F\Mphis[\sigma]F^*=F^*\Mphis[\sigma]F \qand F^*\Mphis[\sigma]
  F^* =F\Mphis[\sigma] F.\] Using this and setting $t=\log(-s)$
  completes the proof.
\end{proof}

\begin{thm}
  \label{thm:uag-FPS}
  The unitary automorphism group of $\FPS$ is generated by
  \[\{ \Mphis ,M_\lambda, F, e^{i\psi}I : s,\lambda,\psi\in\bR\}.\]
\end{thm}
\begin{proof}
  Let $U\in\U(\FPS)$ and let $\Sigmaa$, $\Sigmav$ and $\Sigma_s$ denote the
  ``great circles''
  \[
  \Sigmaa=\anest\cup\anest^\perp,\quad
  \Sigmav=\vnest\cup\vnest^\perp,\quad
  \Sigma_s=\N_s\cup\N_s^\perp\text{ for $s\in\bR$.}
  \]
  The map $\Ad U$ preserves orthogonality and the order structure on $\FPS$, so it must
  permute these great circles.
  If $(\Ad U)\Sigmaa=\Sigmav$,
  then $\Ad FU$ fixes $\Sigmaa$; if $(\Ad U)\Sigmaa=\Sigma_s$ for some $s\in \bR$, then
  $M_{\phis[-s]}U$ fixes $\Sigmaa$. So we may assume that $(\Ad U)\Sigmaa=\Sigmaa$.

  If $(\Ad U)\Sigmav\neq \Sigmav$, then $(\Ad U)\Sigmav=\Sigma_s$ for some $s\neq 0$.  If $s>0$
  then by~(\ref{eq:fourier}), $(\Ad F)\Sigma_s=\Sigma_{-1/s}$
  and $(\Ad F)\Sigmaa=\Sigmav$, so  $U'=F^*M_{\phis[1/s]}FU$ satisfies
  $(\Ad U')\Sigmaa=\Sigmaa$ and
  $(\Ad U')\Sigmav=\Sigmav$.
  So we may assume that $(\Ad U)\Sigmav=\Sigmav$ and $(\Ad U)\Sigmaa=\Sigmaa$.

  There are now four cases to consider:
  \begin{enumerate}
    \item{$(\Ad U)\anest=\anest$, $(\Ad U)\vnest=\vnest$;}\smallskip
    \item{$(\Ad U)\anest=\anest^{\perp}$, $(\Ad U)\vnest=\vnest^{\perp}$;}\smallskip
    \item{$(\Ad U)\anest=\anest^{\perp}$, $(\Ad U)\vnest=\vnest$;}\smallskip
    \item{$(\Ad U)\anest=\anest$, $(\Ad U)\vnest=\vnest^{\perp}$}.
  \end{enumerate}
  Replacing $U$ with $F^2U$ interchanges
  cases~(i) and~(ii) and also interchanges cases~(iii) and~(iv), so it suffices to
  consider cases~(i) and~(iii) only.

  Suppose that case~(iii) holds.
  We claim that $(\Ad U)\N_1=\N_{-s}^\perp$ for some $s>0$.
  To see this, let $\nests$ be the set of nests
  \[\nests=\{ \Nv,\Nv^\perp\}\cup\{ \N_s,\N_s^\perp: s\in\bR\}\]
  so that $\FPS$ is the union of all nests in~$\nests$.  Since $U$ is
  unitary, it maps nests onto nests and so induces a bijection of
  $\nests$.

  Let $\Sigmae$ be the ``equator'' of $\FPS$,
  \[ \Sigmae=\{\Ltp, \Ltm\}\cup \{ \Mphis \Ht, \Mphis\Htb : s\in\bR\}.\]
  Here $\Htb$ is the set of complex conjugates of functions in $\Ht$,
  which is equal to $\Ht^\perp$~\cite{Nik02:vol1}.  The set $\Sigmae$
  contains exactly one subspace from each nest in $\nests$, so the
  action $(\Ad U):\Sigmae\to (\Ad U)\Sigmae$, $K\mapsto (\Ad U)K$ of $\Ad
  U$ on $\Sigmae$ determines the action of $\Ad U$ on $\nests$.
  Moreover, $\Ad U$ is a homeomorphism between $\Sigmae$ and $(\Ad
  U)\Sigmae$, and $\Sigmae$ is itself homeomorphic to the circle~$\bT$.
  Let us give $\nests$ the topology induced by the topology
  on~$\Sigmae$. The bijective action of $\Ad U$ on $\nests$ is then a
  homeomorphism.

  It follows that the closed connected set
  \[
  [\anest,\vnest] = \bigcup_{s\geq 0}\N_s\cup \vnest
  \]
  must be mapped by~$\Ad U$ onto
  \begin{gather*}
    \text{either}\quad
    [\anest^\perp,\vnest]=\bigcup_{s\geq 0}\N_s^\perp\cup
    \vnest^\perp\cup
    \bigcup_{s\in \bR}\N_s\cup \vnest\\[6pt]
    \text{or}\quad
    [\vnest,\anest^\perp] =\vnest\cup \bigcup_{s\leq 0}\N_s^\perp.
  \end{gather*}
  If $(\Ad U)[\anest,\vnest]=[\anest^\perp,\vnest]$ then
  there is some $s>0$ such that $(\Ad U)\N_s=\vnest^\perp$.
  Since~$(\Ad U)\Nv=\Nv$ and $U$ is unitary,
  \[ \Nv^\perp=((\Ad U^*)\Nv)^\perp=(\Ad U^*)\Nv^\perp = \N_s,\]
  which is impossible by the F.~\&~M.~Riesz theorem.

  So $(\Ad U)[\anest,\vnest]=[\vnest,\anest^\perp]$ and so $(\Ad
  U)\N_1=\N_{-s}^\perp$ for some $s>0$.  Since $(\Ad U)\Nv=\Nv$, it
  follows from~\cite{davidson:nest}, Chapter~17 that there
  exist a unimodular function $\alpha\in \Li$ and
  an order-preserving almost everywhere differentiable bijection $g:
  \bR\to\bR$ such that $U=M_\alpha C_g$ where $C_g$ is the unitary
  composition operator corresponding to $g$. Thus \[U
  M_{\phis[1]}\Ht=M_\alpha C_g M_{\phis[1]}\Ht=M_{\phis[-s]}M_\lambda
  \Htb\] for some $\lambda\in \bR$. Moreover, $(\Ad
  U)\anest=\anest^\perp$, so $U\Ht=M_\mu \Htb$ for some $\mu\in \bR$.
  Since $C_gM_f=M_{f\circ g}C_g$ for $f\in\Li$,
  \begin{equation*}
    M_{\phis[1]\circ g}M_\alpha C_g \Ht = M_{\phis[1] \circ g}M_\mu \Htb
    =M_{\phis[-s]}M_\lambda \Htb.
  \end{equation*}
  Taking orthogonal complements, we see that
  \begin{equation*}
    \label{eq:Htb-reducing}
    \Mphis M_{-\lambda} M_{\phis[1]\circ g}M_\mu \Ht=u\Ht=\Ht,
  \end{equation*}
  where $u:\bR\to\bC$ is the unimodular function
  \[x\mapsto\exp i\big( -\tfrac12(g(x)^2+sx^2)+(\mu-\lambda)x
  \big).\] So $u$ must be constant almost everywhere. But
  $g(x)\to \infty$ as $x\to \infty$, so this is impossible.

  So we are reduced to case~(i): $\Ad U$ fixes both the analytic nest and
  the Volterra nest, and so is a unitary automorphism of
  $\Alg(\vnest\cup\anest)$, the Fourier binest algebra.
  By~\cite{pow-kat:FB}, Lemma~4.1, $U=e^{i\psi}M_\lambda
  D_\mu V_t$ for some $\psi,\lambda,\mu,t\in\bR$. Now apply
  Lemma~\ref{lem:Lang-Vt}.
\end{proof}

\begin{rem}
  It can be shown that, modulo scalars, this automorphism group is
  isomorphic to the semidirect product $\bR^2\rtimes SL_2(\bR)$. The
  isomorphism is implemented by the map sending
  \[
  \begin{bmatrix}
    1 \\
    \lambda & 1 & 0\\
    \mu & 0 & 1
  \end{bmatrix},\quad
  \begin{bmatrix}
    1 \\
     & 1 & -s\\
     & 0 & 1
  \end{bmatrix}\qand
  \begin{bmatrix}
    1 \\
     & 0 & -1\\
     & 1 & 0
  \end{bmatrix}
  \]
  to $\Ad(M_\lambda D_\mu)$, $\Ad(\Mphis)$ and $\Ad(F)$ respectively.
  We refer the reader to~\cite{levene-thesis} for the details.  

  It is perhaps surprising that $\{\Ad(U):U\in\U(\FPS)\}$
  has such a simple description. The authors do not know if the same
  can be said for $\U(\FPS)$ itself.
\end{rem}

\subsection{The hyperbolic sphere and the extended hyperbolic lattice}
\label{sec:5.2}

Recall the definitions of the hyperbolic sphere $\Shyp$ and
the extended hyperbolic lattice $\Mhyp\supseteq\Shyp$ from
Example~\ref{eg:hyp-sphere}. For $\lambda,\mu\in\bR$, let
$M_{\lambda,\mu}=M_{e^{i(\lambda x+\mu x^{-1})}}$
and for $(\theta,s)\in\bT\times\bR$ let
\[U_{\theta,s}=M_{|x|^{is}u_\theta(x)}\quad\text{where as before,}\quad
u_\theta(x)=\chi_{[0,\infty)} + \theta\chi_{(-\infty,0)}.\]
A typical projection
in $\M(\S_{(x,\log|x|,-x^{-1})})$ is \[[U_{1,s}M_{\lambda,\mu}\Ht] \q{where}
(s,\lambda,\mu)\in\bR^3.\]
If $u_{\theta,s}=\chi_{[0,\infty)} + \theta
e^{s\pi}\chi_{(-\infty,0)}$ then it is shown
in~\cite{pow-kat:hyp} that
$U_{\theta,s}\Ht = M_{u_{\theta,s}}\Ht$ for $(\theta,s)\in\bT\times\bR$.
We further define operators
\begin{equation*}
  J_1 f(x)=x^{-1}f(-x^{-1})\qand
  J_2f(x)=f(-x), \quad f\in\Lt;
\end{equation*}
these are the unitary composition operators corresponding to
the symmetries $x\mapsto -x^{-1}$ and $x\mapsto -x$, respectively.
The linear span of the set of functions $z\mapsto (z-\xi)^{-1}$ for $\Im
\xi < 0$ (or $\Im \xi > 0$) is dense in $\Ht$ (or in $\Htb$,
respectively). Applying $J_1$ and $J_2$ to these sets reveals that
$J_1 \Ht=\Ht$ and $J_2 \Ht=\Htb$.  It is easy to see that
all of these operators are unitary automorphisms of $\Mhyp$, and if we
fix $\theta=1$ then we obtain unitary automorphisms of $\Shyp$.  We
will show that in each case, these operators generate the whole
unitary automorphism group.

\begin{lem}
  \label{lem:gamma}
  Let $\gamma$ be a conformal automorphism of the upper half plane
  $\UHP$. For each nonzero $s\in\bR$,
  the
  subspace $M_{u_{1,s}\circ \gamma} \Ht$ has zero intersection with
  each subspace in $\M(\S_{(x,\log|x|,-x^{-1})})$
  unless $\gamma$ is either of the form $\gamma(x)=ax$
  for some $a>0$ or $\gamma(x)=-bx^{-1}$ for some $b>0$.
\end{lem}
\begin{proof}
  Suppose that $M_{u_{1,s}\circ \gamma} \Ht \cap
  U_{1,\sigma}M_{\lambda,\mu}\Ht\neq \{0\}$. Let \[f=(u_{1,s}\circ
  \gamma) g = u_{1,\sigma}e^{i(\lambda x+\mu x^{-1})}h\] be a
  nonzero function in this
  intersection, where $g,h\in \Ht$. Multiplying this equation  by
  $e^{-i\lambda x}$ if $\lambda < 0$ and by $e^{-i\mu x^{-1}}$ if $\mu >0$
  and writing $\alpha=(u_{1,s}\circ
  \gamma)/u_{1,\sigma}$ gives $\alpha\phi = \psi$ for nonzero functions $\phi,\psi\in
  \Ht$.
  Observe that $\alpha$ takes at most four values, since
  \[\alpha(x)=
  \begin{cases}
    1\quad &x\in \gamma^{-1}(\bR_+)\cap \bR_+\\
    e^{s\pi} & x\in \gamma^{-1}(\bR_-)\cap \bR_+\\
    e^{-\sigma\pi} & x\in \gamma^{-1}(\bR_+)\cap \bR_-\\
    e^{(s-\sigma)\pi} & x\in \gamma^{-1}(\bR_-)\cap \bR_-
  \end{cases}.\]
  Applying the F.~\&~M.~Riesz theorem to the function
  $\alpha\phi-\psi=0$ reveals that $\alpha$ must be constant almost
  everywhere; in particular, since $s\ne0$, no three of these
  intersections can have nonzero Lebesgue measure.

  Since $\gamma$ induces a conformal automorphism of the upper half plane,
  either $\gamma(x)=ax+b$ with $a>0$ and $b\in\bR$ or
  $\gamma(x)=a-b(x-c)^{-1}$ with $a,b,c\in\bR$ and $b>0$. Applying the
  condition in the previous paragraph forces $\gamma(x)=ax$ with $a>0$ or
  $\gamma(x)=-b x^{-1}$ with $b>0$.
\end{proof}

\begin{thm}
  \label{thm:autgroups}
  (i) The unitary automorphism group of $\Mhyp$ is equal to the union
  $G\cup GJ_1\cup GJ_2 \cup GJ_1J_2$ where
  \[ G=\{\alpha U_{\theta,s} M_{\lambda,\mu} V_t :
  (\alpha,\theta,s,\lambda,\mu,t)\in \bT^2\times \bR^4\}.\]
  (ii) The unitary automorphism group of $\Shyp$ is equal to the
  union   $G_0\cup G_0J_1\cup G_0J_2 \cup G_0J_1J_2$ where
  \[ G_0=\{\alpha U_{1,s} M_{\lambda,\mu} V_t :
  (\alpha,s,\lambda,\mu,t)\in \bT\times \bR^4\}.\]
\end{thm}
\begin{proof}
  (i) We exploit the order structure of $\Mhyp$, given in
  Proposition~5.2 of~\cite{pow-kat:hyp}.  Suppose that $U$ lies
  in $\Uhyp$, the unitary automorphism group of $\Mhyp$. Let us
  write $\M$ for the set
  \[\M=\M(\S_{(x,\log|x|,-x^{-1})})
  =\{\Ad(U_{1,s}M_{\lambda,\mu})[\Ht]:s,\lambda,\mu\in\bR\}\]
  and let $\bdy$ be the topological boundary of $\M$, which by
  Theorem~\ref{thm:sot-closures} is the set of projections in
  $\Mhyp$ of the form $[L^2(E)]$.  Observe first that $\bdy$ must
  be mapped onto itself by $\Ad U$. This will follow if we can
  show that the set $\Mhyp\setminus\bdy$ may be intrinsically
  described as the union of all non-commutative sublattices of
  $\Mhyp$ which are order-isomorphic to the slice $\L_{1,0}=\{
  \Ad(M_{\lambda,\mu})[\Ht]:\lambda,\mu\in\bR\}$ and whose closure
  contains $\{0,I\}$. Writing
  \[\L_{\theta,s}=\{\Ad(U_{\theta,s}M_{\lambda,\mu})[\Ht]:\lambda,\mu\in\bR\}\]
  for $\theta,s\in\bR$, observe that each of the slices
  $\L_{\theta,s}$, $\L_{\theta,s}^\perp$ has this property. Hence
  this union contains $\Mhyp\setminus\bdy$. On the other hand,
  suppose that $\L$ is such a non-commutative sublattice and that
  $P\in \bdy\cap\L$. Since $\L\cong\L_{1,0}$ there are two
  continuous nests $\N_1,\N_2$ contained in $\L\cup\{0,I\}$ which
  do not commute with one another such that
  $\N_1\cap\N_2=\{0,P,I\}$. Now $P$ is of the form $P=[L^2(E)]$
  and all the non-trivial projections $Q,R$ in $\Mhyp$ which
  satisfy $Q\le [L^2(E)]\le R$ are all contained in $\bdy$ by the
  F.~\&~M.~Riesz theorem, so $\N_1,\N_2\subseteq\bdy$. But all
  projections in $\bdy$ commute, so we obtain a contradiction and
  $\bdy\cap\L$ must be empty.

  Hence $(\Ad U)\bdy=\bdy$. Observe that
  \[\Mhyp\setminus\bdy=\bigcup_{\theta\in\bT}(\Ad U_{\theta,1})\M\cup
  \bigcup_{\theta\in\bT} (\Ad U_{\theta,1})\M^\perp\] and that the
  terms in this union are the components of $\Mhyp\setminus\bdy$,
  which $\Ad U$ must therefore permute. If
  $U\Ht=U_{\theta,s} M_{\lambda,\mu}\Htb$, then we may
  replace $U$ with $J_2U$ to ensure that
  $U\Ht\in\bigcup_{\theta\in\bT}(\Ad U_{\theta,1})\M$, and then
  replacing $U$ by $U_{\theta,s}M_{\lambda,\mu}U$ for suitable
  $\theta,s,\lambda,\mu$ we can arrange that $U\Ht=\Ht$, and so
  also that $(\Ad U)\M=\M$.

  Since $(\Ad U)\bdy=\bdy$, $\Ad U$ maps projections in $(\bdy)''$ to
  projections in $(\bdy)''$ and so induces an automorphism of
  $\Li$. Von Neumann's theorem of~\cite{vonNeumann} shows that this is
  necessarily induced by a Borel isomorphism $\gamma$.  (See also
  Nordgren~\cite{nordgren}).  Since $\Li$ is maximal abelian it
  follows readily that $U=M_\phi C_\gamma$ for some unimodular
  function $\phi$ where $C_\gamma$ is the unitary composition operator
  for $\gamma$.  Moreover, $U$ induces an automorphism of
  $H^\infty(\bR)$, since if $h\in H^\infty(\bR)$ then $M_h \Ht
  \subseteq \Ht$ and so
  \begin{align*} \Ht = U\Ht \supseteq U M_h \Ht & = M_\phi
    C_\gamma M_h \Ht \\& = M_{h\circ \gamma} M_\phi C_\gamma \Ht \\&=
    M_{h\circ \gamma} U\Ht \\& =M_{h\circ \gamma} \Ht.\end{align*} It
  follows from the Beurling-Lax theorem~\cite{lax59:translation}
  that $h\circ\gamma\in H^\infty(\bR)$. The same argument applied
  to $U^*$ shows that the map $h\mapsto h\circ\gamma$ is indeed
  an automorphism of $H^\infty(\bR)$. Hence $\gamma$ induces a
  conformal automorphism of the upper half plane.

  Fix $s\neq 0$. Since $(\Ad U)\M=\M$, the subspace
  \[ U (U_{1,s}\Ht) = U M_{u_{1,s}}\Ht = M_{u_{1,s}\circ \gamma} U\Ht
  = M_{u_{1,s}\circ \gamma} \Ht\] must lie in $\M$. By
  Lemma~\ref{lem:gamma}, either $\gamma(x)= e^t x$ or
  $\gamma(x)=-(e^t x)^{-1}$ for some $t\in \bR$.
  Hence $C_\gamma\in \{
  V_t, J_1V_t\}$.  Multiplying by $C_\gamma^*$ reduces to the case
  $U=M_\phi$. Since $C_\gamma \Ht=\Ht$, we have $M_\phi \Ht=\Ht$ and so $\phi$ is
  constant almost everywhere.

  (ii) Following the first paragraph of the proof of~(i) with
  $\Shyp$ in place of $\Mhyp$, we may assume that $U\in
  \U(\Shyp)$ satisfies $(\Ad U)\bdy=\bdy$. Now
  $\Shyp\setminus\bdy=\M\cup\M^\perp$ has two components, which
  $\Ad U$ must permute. We may again assume that $(\Ad
  U)\M=\M$ by multiplying by $J_2$ if necessary. The remainder of
  the proof proceeds as above.
\end{proof}

From Table~\ref{table:comm},
\begin{table}
  \caption{Commutation relations for $\Uhyp$}
  \label{table:comm}
  \centering
  \begin{tabular}{rr|cccccc}
    & & \multicolumn{5}{c}{$Y$} \\%\hline
    \multirow{7}*{$X$} & $XY$
    & $U_{\theta,s}$           & $M_{\lambda,\mu}$ & $V_t$ & $J_1$
    & $J_2$ %& $J_0$
    \\\hline
    &$U_{\theta',s'}$ & $U_{\theta\theta',s+s'}$ & \\%\hline
    &$M_{\lambda',\mu'}$ & commute                  &
    $M_{\lambda+\lambda',\mu+\mu'}$ \\%\hline
    &$V_{t'}$ & $e^{ist'}U_{\theta ,s}V_{t'}$         &
    $M_{e^{t'}\lambda,e^{-t'}\mu}V_{t'}$   & $V_{t+t'}$ \\%\hline
    &$J_1$ & $U_{\overline\theta,-s}J_1$       & $M_{-\mu,-\lambda}J_1$          &
    $V_{-t}J_1$ & $I$ & % & $J_2$
    \\%\hline
    &$J_2$ & $U_{\overline{\theta},s}J_2$      & $M_{-\lambda,-\mu}J_2$          &
    commute & commute & $I$ %& $J_1$
    \\%\hline
  \end{tabular}
\end{table}
we see that $\Ad(\Uhyp)$ is isomorphic to the double semidirect
product
\begin{align*}\big( (\bT\times \bR^3) \rtimes_\alpha \bR\big)&\rtimes_\beta
(\bZ/2\bZ)^2,\\
\alpha(t)(\theta,s,\lambda,\mu)&=(\theta,s,e^t\lambda,e^{-t}\mu),\\
\beta(1,0)(\theta,s,\lambda,\mu,t)&=(\overline{\theta},-s,-\mu,-\lambda,-t),\\
\beta(0,1)(\theta,s,\lambda,\mu,t)&=(\overline{\theta},s,-\lambda,-\mu,t).\end{align*}
The map sending
\[
\begin{bmatrix}
  1        &   &   &         & \\
  a        & 1 & 0 &         & \\
  s        & 0 & 1 &         & \\
  \lambda  &   &   &  e^{t} &0\\
  \mu      &   &   &    0    & e^{-t}
\end{bmatrix},\!\quad\!
\begin{bmatrix}
  -1 &   &   &   &   \\
     & 1 & 0 &   &   \\
     & 0 & 1 &   &   \\
     &   &   & 0 & 1 \\
     &   &   & 1 & 0 \\
\end{bmatrix}\!\qand\!
\begin{bmatrix}
  -1 &   &   &   &   \\
     & 1 & 0 &   &   \\
     & 0 & -1 &   &   \\
     &   &    & 1 & 0 \\
     &   &    & 0 & 1 \\
\end{bmatrix}
\]
to $\Ad(U_{e^{ia},s}M_{\lambda,\mu}V_t)$, $\Ad(J_1)$ and
$\Ad(J_2)$ respectively for $a,s,\lambda,\mu,t\in \bR$ is a
homomorphism onto $\Ad(\Uhyp)$ with kernel $2\pi\bZ$.

\subsection{Polynomial $2$-spheres}
\label{sec:5.3}

Finally, we consider isomorphisms between
$2$-spheres of the form $\Sigma_{m,n}=\Sigma(\M_{m,n})$ where
\[ \M_{m,n}=
\M(\S_{(x^m,x^n)})\q{for}m,n\in\bZ\setminus\{0\}.\]
We write $\bdy_{m,n}$ for the topological boundary of $\M_{m,n}$.
This can be easily computed using the results of Section~\ref{sec:4}:

\begin{lem}
  \label{lem:bdys}
  Let $m>n$ be nonzero integers and let $\alpha$ be the residue class
  of $(m,n)$ in $(\bZ/2\bZ)^2$, which we identify with
  $\{0,1\}\times\{0,1\}$. Then $\bdy_{m,n}=\bdy_\alpha$ depends only
  on $\alpha$. In fact $\bdy_{(0,1)}=\Sigmav$, the Volterra circle,
  \[\bdy_{(1,1)}=\{P,P^\perp:P=[\Ltwo{E}],\ E=[-a,a],\ a\in
  [0,\infty]\}\] and if $\beta=(1,1)-\alpha$ then
  \[\bdy_\beta =
  \{PM_{\chi_{(0,\infty)}}+P^\perp M_{\chi_{(-\infty,0)}}:P\in
  \bdy_\alpha\}.
  \]
\end{lem}

We now examine the order structure on $\Sigma_{m,n}$, which is rather
simple in many cases.

\begin{lem}
  \label{lem:noblaschke}
  Suppose that $u(x) = \exp(i(\gamma x^m + \delta x^n))$ is an inner
  function, where $n,m$ are distinct integers and $\gamma,\delta\in\bR$. If
  $\gamma\ne0$ then $m\in \{-1,0,1\}$, and if $\delta \ne 0$ then
  $n\in \{-1,0,1\}$.
\end{lem}
\begin{proof}
  Given $h(x)\in H^\infty(\bR)$,
  there is a unique function, $\eta(z)\in H^\infty(\UHP)$
  whose nontangential limit boundary value function
  $\bfn{\eta}\in H^\infty(\bR)$ is equal to $h$ almost
  everywhere. Suppose that $\kappa(z)$ is analytic on an open disc $U$
  with $U \cap \bR = (a,b)$ for some $a<b$, and that 
  $h$ agrees with $\bfn{\kappa}$ almost everywhere on
  $(a,b)$. Then $\eta(z) = \kappa(z)$ on $U \cap \UHP $. Indeed $\eta$ and $\kappa$ both
  have restrictions in $H^\infty (U \cap \UHP)$ and their boundary
  functions agree on a set of positive measure,  so the conclusion
  follows from the F.~\&~M.~Riesz theorem together with the Riemann
  mapping theorem.  (See also Fisher~\cite{fisher}.)

  We apply this principle to $h(x) = \exp(i(\gamma x^m + \delta
  x^n))$. If $h$ is inner and $\eta\in H^\infty(\UHP)$ with
  $\bfn{\eta}=h$, consider the analytic function
  $\kappa:\bC\setminus\{0\}\to\bC$, $z\mapsto\exp(i(\gamma z^m +
  \delta z^n))$. Since $\kappa$ is analytic on each of the open discs
  $U$ which meet the right half line and the union of these discs
  contains~$\UHP$, we conclude that $\eta=\kappa$ on~$\UHP$.  However,
  it is routine to check that this function is bounded in the upper
  half plane only under the stated conditions.
\end{proof}

\begin{lem}
  \label{lem:polysphere-order}
  Let $m,n$ be distinct integers in $\bZ\setminus\{-1,0,1\}$.  If $P$
  and $Q$ are distinct projections in $\Sigma_{m,n}$ with $0\ne P\leq Q\ne
  I$ then they must lie in $\bdy_{m,n}$.
\end{lem}
\begin{proof}
  Suppose that $P\not\in\bdy_{m,n}$. By the F.~\&~M.~Riesz theorem,
  $Q\not\in\bdy_{m,n}$. Suppose without loss of generality that
  $P\in\M_{m,n}$; if this is not the case, then apply the unitary
  automorphism $J_2$
  induced by $x\mapsto -x$, which maps $\M_{m,n}^\perp$ onto
  $\M_{m,n}$, to make it so. If
  $Q\in\M_{m,n}^\perp$ then there exist
  $\alpha,\beta,\lambda,\mu\in\bR$ such that
  \[ e^{i(\lambda x^m+\mu x^n)}\Ht \subseteq
  e^{i(\alpha x^m+\beta x^n)}\Htb,\]
  where these subspaces are the ranges of $P$ and $Q$,
  respectively. Let $u(x)=\exp{i(\gamma x^m+ \delta x^n)}$ where
  $\gamma=\lambda-\alpha$ and $\delta=\mu-\beta$, so that
  $u\Ht\subseteq\Htb$. Taking complex conjugates yields $\overline
  u\Htb \subseteq \Ht$, so $\Htb\subseteq u\Ht$ and thus
  $u\Ht=\Htb$. So
  \[\Ht=(\Htb)^\perp=(u\Ht)^\perp = u\Htb = u^2\Ht.\]
  It is
  well-known that a unimodular function which preserves $\Ht$ must be
  constant almost everywhere, so $\gamma=\delta=0$, which would imply that
  $\Ht=\Htb$, an obvious contradiction.

  So $Q\in \M_{m,n}$, say \[
  P=[e^{i(\lambda x^m+\mu x^n)}\Ht],\quad Q= [e^{i(\alpha x^m+\beta
    x^n)}\Ht].\]
  Now $u(x)=\exp{i(\gamma x^m+ \delta x^n)}$ leaves $\Ht$ invariant,
  so is an inner function. Hence $\gamma=\delta=0$ by
  Lemma~\ref{lem:noblaschke} and $P=Q$, a contradiction.
\end{proof}

The conclusion of Lemma~~\ref{lem:polysphere-order} holds nontrivially
precisely when $m, n$ are not both even and it follows that the
boundary $\bdy_{m,n}$ is a unitary invariant for $\Sigma_{m,n}$ and
also for the balls $\M_{m,n}$.  Furthermore, Lemma~\ref{lem:bdys}
shows that $\bdy_{m,n}$ generates $\Li$ precisely when $m, n$ are not
both odd.  When both these conditions prevail we can classify the
spheres and balls by an argument similar to that of
Theorem~\ref{thm:autgroups}. In fact we expect that a somewhat deeper
analysis will show that in general the (unordered) set
\[\big\{\{m,n\},\{-m,-n\}\big\}\] is a complete unitary invariant.
\medskip

We shall need the following elementary lemma.

\begin{lem}
  \label{lem:halfplane}
  Let $m,n,p,q$ be nonzero integers with $m,n\ge1$ such that $p\ne
  q$ and $m\ne n$. 
  If $\gamma:\bR\to\bR$ induces a conformal automorphism of
  the upper half plane and $\alpha,\beta$ are real constants
  such that
  \[ x^p-x^q=\alpha \gamma(x)^m+\beta\gamma(x)^n + 2\pi N(x)\q{almost
    everywhere,}\] where $N:\bR\to\bZ$ then
  \[\big\{\{m,n\},\{-m,-n\}\big\}=\big\{\{p,q\},\{-p,-q\}\big\}.\]
\end{lem}
\begin{proof}
  Without loss of generality, suppose that $p>q$ and $m>n$.  Either
  $\gamma(x)=ax+b$ where $a>0$ and $b\in\bR$, or
  $\gamma(x)=a-b(x-c)^{-1}$ for $a,c\in\bR$ and $b>0$.  Suppose first
  that $\gamma(x)=ax+b$; without loss of generality, we may take
  $a=1$. Since $N$ is then continuous and so constant on $(0,\infty)$
  it follows that $p,q\ge1$ and so $N$ is constant on~$\bR$. The
  equation
  \[ x^p-x^q = \alpha(x+b)^m + \beta(x+b)^n +2\pi N\]
  holds almost everywhere. Considering the
  coefficient of $x^p$ gives $\alpha=1$ and $m=p$, so we suppose that
  $n\ne q$.  Differentiating
  gives
  \[ px^{p-1}-qx^{q-1} = p(x+b)^{p-1} + \beta n(x+b)^{n-1}. \] If
  $q,n>1$ then we set $x=-b$ to deduce that $b$ is algebraic, and set
  $x$ equal to any other algebraic number to see that $\beta$ is also
  algebraic. Simple arguments show that the same holds if $q>1$ and
  $n\ne 1$ or if $q=1$ and $n>1$. Now equate the constant terms in the
  original expression:
  \[2\pi N = -(b^p + \beta b^n).\] Since the right hand side is
  algebraic, $N=0$. By counting repeated roots, it now follows that
  $n=q$.

  If on the other hand $\gamma(x)=a-b(x-c)^{-1}$ then $N$ is
  continuous and so constant on $\bR\setminus\{0,c\}$, so $p,q\le-1$
  and $N(x)=0$ almost everywhere for $x>\max\{0,c\}$ and for
  $x<\min\{0,c\}$.  Since the left hand side is locally unbounded
  only at $x=0$ and has limit $0$ as $x\to\pm\infty$, we must have
  $c=0$ and $N(x)=0$ almost everywhere.  It only remains to consider
  the order of growth and decay at $0$ and $\pm\infty$ to see that
  $p=-n$ and $q=-m$.
\end{proof}

\begin{thm}
  Let $p, q\in \bZ\setminus\{-1,0,1\}$ with $p\ne q$ and let
  $m, n > 1$ be integers with $m\not\equiv n\mod 2$. 
  The spheres  $\Sigma_{m,n}$ and $\Sigma_{p,q}$ are unitarily
  equivalent if and only if
  \[\big\{\{m,n\},\{-m,-n\}\big\}
  = \big\{\{p,q\},\{-p,-q\}\big\}.
  \]
\end{thm}
\begin{proof}
  First observe that the spheres are unitarily equivalent if these
  sets are equal, since the composition operator $J_1$ corresponding to
  the map $x\mapsto -x^{-1}$ satisfies $(\Ad
  J_1)\Sigma_{m,n}=\Sigma_{-m,-n}$.

  Let $U\in\Unit(\Lt)$ with $(\Ad U)\Sigma_{m,n}=\Sigma_{p,q}$.
  Consider the subspace $U\Ht\in\Sigma_{p,q}$. Since $\Ad U$ preserves
  the order structure, it must map $\bdy_{m,n}$ onto $\bdy_{p,q}$ by
  Lemmas~\ref{lem:bdys} and~\ref{lem:polysphere-order}. By composition
  with $x\mapsto -x$ if necessary, we may assume that $U\Ht\in
  \M_{p,q}$ and then translating by the ``obvious'' inner
  automorphisms of $\M_{p,q}$, that $U\Ht=\Ht$.

  Let $\sigma$ be the residue class of
  $(m,n)$ in $(\bZ/2\bZ)^2$, which by assumption is in
  $\{(0,1),(1,0)\}$. Observe that by Lemma~\ref{lem:bdys}, the
  von Neumann algebra generated by $\bdy_{m,n}=\bdy_\sigma$ is the
  multiplication algebra $\Li$, and the only possibilities for
  the algebra $\A=(\bdy_{p,q})''$ are
  \[\A=\Li\qor\A=\{M_f:f\in \Li,\ f(x)=f(-x)\}.\]
  The latter algebra has uniform multiplicity~$2$.  Since $\Ad U$
  sends projections in $(\bdy_{m,n})''$ to projections in
  $(\bdy_{p,q})''$, it induces an isomorphism between $\Li$ and $\A$.
  Spatial isomorphisms preserve multiplicity, so in fact $\A=\Li$.

  Now $\Ad U$ is an isomorphism $\Li\to\Li$ and it follows that
  $U=M_\phi C_\gamma$ where $\phi\in\Li$ is unimodular and $C_\gamma$
  is a unitary composition operator with symbol~$\gamma$, a Borel
  isomorphism.  Exactly as in the proof of
  Theorem~\ref{thm:autgroups}, $\gamma$ induces a conformal
  automorphism of the upper half plane. Since $\Ad U$ is a
  homeomorphism with \[U\Ht=\Ht\qand U\bdy_{m,n}=\bdy_{p,q},\] it maps
  the two components $\M_{m,n}$ and $\M_{m,n}^\perp$ of
  $\Sigma_{m,n}\setminus\bdy_{m,n}$ to $\M_{p,q}$ and $\M_{p,q}^\perp$
  respectively. In particular, there exist real $\alpha,\beta$ such
  that  \begin{align*}
    e^{i(x^p-x^q)}\Ht
    &= Ue^{i(\alpha x^m+\beta x^n)}\Ht\\
    &= e^{i(\alpha \gamma(x)^m+\beta \gamma(x)^n)}U\Ht\\
    &=e^{i(\alpha \gamma(x)^m+\beta \gamma(x)^n)}\Ht.
  \end{align*}
  Hence the hypotheses of
  Lemma~\ref{lem:halfplane} are satisfied and the result follows.
\end{proof}

\bibliographystyle{abbrv}
\def\lfhook#1{\setbox0=\hbox{#1}{\ooalign{\hidewidth
  \lower1.5ex\hbox{'}\hidewidth\crcr\unhbox0}}}

\end{document}